\documentclass[numbib]{imamat}
\usepackage{mathtools}
\usepackage{amsmath}	% MATH PACKAGE INCLUDING THM DISPLAYING
\usepackage{amssymb}	% MATH PACKAGE INCLUDING MATH SYMBOLS AND FONTS
\usepackage{float}

\usepackage{amsthm}
\usepackage{graphicx,graphics,color}	% FOR FIGURES
\usepackage{enumerate}
\usepackage{cancel}%strike out text with \sout
\usepackage{bm}

\usepackage{epsfig}
\usepackage{subfigure}
\usepackage{xcolor}
\newcommand{\red}{\color[rgb]{1,0,0}}
\newcommand{\green}{\color[rgb]{0,1,0}}
\newcommand{\blue}{\color[rgb]{0,0,1}}

\usepackage{setspace}
\usepackage{dcolumn}% Align table columns on decimal point
\usepackage{bm}% bold math
\usepackage{times}

\usepackage[titletoc,title]{appendix}
\usepackage[titletoc]{appendix}

\usepackage[titletoc,title]{appendix}

\theoremstyle{proof} \newtheorem{assume}{Assumption}[section]
\theoremstyle{proof}

\DeclareFontFamily{U}{msb}{}
\DeclareFontShape{U}{msb}{m}{n}{ <5> <6> <7> <8> <9> gen * msbm
  <10> <10.95> <12> <14.4> <17.28> <20.74> <24.88> msbm10}{} 
\DeclareSymbolFont{AMSb}{U}{msb}{m}{n}
\DeclareMathSymbol{\N}{\mathalpha}{AMSb}{"4E}
\DeclareMathSymbol{\R}{\mathalpha}{AMSb}{"52}
\DeclareMathSymbol{\Z}{\mathalpha}{AMSb}{"5A}

\newcommand{\tlamb}{{\widetilde{\lambda}}}

\newcommand{\bu}{\bm{u}}
\newcommand{\tu}{{\widetilde{\bm{u}}}}
\newcommand{\bx}{\bm{x}}
\newcommand{\bC}{\bm{C}}
\newcommand{\bX}{\bm{X}}

\DeclareRobustCommand{\E}[1]{{\mathbb{E}}\left[{#1}\right]}
\DeclareRobustCommand{\wE}[1]{{\widehat{\mathbb{E}}}\left[{#1}\right]}

\begin{document}

\title{Ensemble-based estimates of eigenvector error for empirical covariance matrices}

\shorttitle{ENSEMBLE-BASED ESTIMATES OF EIGENVECTOR ERROR} %%%for recto running head
\shortauthorlist{D. TAYLOR, J. G. RESTREPO, AND F. G. MEYER} %%% for verso running head

\author{%%%% First author details
\name{Dane Taylor}
\address{
Department of Mathematics, University at Buffalo, State University of New York, Buffalo, NY 14260, USA; Statistical and Applied Mathematical Sciences Institute, Research Triangle Park, NC 27709, USA; Department of Mathematics, University of North Carolina, Chapel Hill, NC 27599, USA\email{danet@buffalo.edu}}
%%%%%%% Second author details
\name{Juan G. Restrepo}
\address{Department of Applied Mathematics, University of Colorado, Boulder, CO, 80309, USA
\email{juanga@colorado.edu}}
%%%%%%%
\and
%%%%%%% Third author details
\name{Fran\c cois G.  Meyer}
\address{Department of Electrical, Computer and Energy Engineering, University of Colorado, Boulder, CO, 80309, USA}
\email{fmeyer@colorado.edu}}

\maketitle

\begin{abstract}
{
Covariance matrices are fundamental to the analysis and forecast of economic, physical and biological systems. Although the eigenvalues $\{\lambda_i\}$ and eigenvectors $\{\bm{u}_i\}$ of a covariance matrix are central to such endeavors, in practice one must inevitably approximate the covariance matrix based on data with finite sample size $n$ to obtain empirical eigenvalues $\{\tilde{\lambda}_i\}$ and eigenvectors $\{\tilde{\bm{u}}_i\}$, and therefore understanding the error so introduced is of central importance.  We analyze eigenvector error $\|\bm{u}_i - \tilde{\bm{u}}_i \|^2$ while leveraging the assumption that the true covariance matrix having size $p$ is drawn from a matrix ensemble with known spectral properties---particularly, we assume the distribution of population eigenvalues weakly converges as $p\to\infty$ to a spectral density $\rho(\lambda)$ and that the spacing between population eigenvalues is similar to that for the Gaussian orthogonal ensemble. Our approach complements previous analyses of eigenvector error that require the full set of eigenvalues to be known, which can be computationally infeasible when $p$ is large. To provide a scalable approach for uncertainty quantification of eigenvector error, we consider a fixed eigenvalue $\lambda$ and approximate the distribution of the expected square error $r= \mathbb{E}\left[\| \bm{u}_i - \tilde{\bm{u}}_i \|^2\right]$ across the matrix ensemble for all $\bm{u}_i$ associated with $\lambda_i=\lambda$. We find, for example, that for sufficiently large matrix size $p$ and sample size $n> p$, the probability density of $r$ scales as $1/nr^2$. This power-law scaling implies that  eigenvector error is extremely heterogeneous---even if $r$ is very small for most eigenvectors, it can be  large for others with non-negligible probability. We support this and further results with numerical experiments.
}%
{covariance matrix, empirical eigenvector, Wigner surmise, Wishart distribution, Graphical model.}
\\
Math Subject Classification: 37N99, 62B10, 94A17
%37N99, 62B10, 94A17
%% 37N99: Dynamical systems and ergodic theory -> Applications
%% 62B10: Statistics -> Information-theoretic topics
%% 94A17: Information and communications, circuits -> Measures of information, entropy
\end{abstract}

%________________________________________________________________
\section{Introduction}\label{sec:intro}
%________________________________________________________________

The spectral properties of covariance matrices are a central topic in mathematics, probability and statistics (\cite{Golub2012,anderson03,Mehta1991,Hastie2009}) and provide a cornerstone to applications in physics, biology, economics and social science (\cite{Elton2009,Mantegna2000,Delvenne2010,Volkov2009,Bassett2011,Gatti2010,Weigt2009}).  The estimation of eigenvectors of a sample covariance matrix remains a fundamental tool for these and numerous other application domains. Sample covariance matrices can be computed locally if the dataset lies along a manifold, or globally if the data is organized along a linear subspace (\cite{Hastie2009}). Often, the practitioner has access to a generative stochastic model for the covariance matrix that can be derived from first principles or domain knowledge, and he/she needs to estimate the accuracy of eigenvectors calculated from a sample covariance matrix.

{We consider the ``classical'' (large sample, $n>p$) framework where one has access to $n$ measurements of a $p$-dimensional vector $\bx$ with covariance $\bC$, which are encoded as the columns of a matrix $\bX$ of size $p\times n$.
% is the size-$(p\times n)$ matrix having the $\bx^{(i)}$ as columns
%
%which are  drawn as i.i.d. samples from a multivariate normal distribution $\mathcal{N}_p(0,\bC)$ with zero mean and population covariance matrix $\bC$. The maximum-likelihood estimate  for $\bC$ is well known to be given 
%
The sample covariance matrix $\tilde{\bC}= n^{-1} (\bX -\mathbb{E}[\bX ])(\bX-\mathbb{E}[\bX])^T$ is an unbiased estimator to the population covariance matrix $\bC$, and
%, where $\bX$ is the size-$(p\times n)$ matrix having the $\bx^{(i)}$ as columns (\cite{anderson03}).
%
the main motivation for our work is to estimate how well the eigenvectors $\{ \tu_i \}$ of $\tilde{\bC}$ approximate the eigenvectors $\{\bu_i\}$ of $\bC$ in the limit when both $p$ and $n$ are large.
If we further assume that $\tilde{\bC}$ is distributed according to a Wishart distribution $W(\bC,n)$ centered at $\bC$ (which occurs, for example, when $\bx$ follows a multivariate normal distribution), then
for fixed $p$ and $n \rightarrow \infty$, the expected error between a sample eigenvector $\tu_i$
% for a Wishart-distributed matrix $\tilde{\bC}$ and the corresponding population eigenvector 
and $\bu_i$ for $\bC$ for $i\in\{1,\dots,p\}$ 
is given asymptotically by (\cite{anderson03}, Theorem 13.5.1) 
%\todo{Francois, WHAT THEOREM or PAGE??????}
\begin{equation}
  \E{n \|\bu_i - \tu_i\|^2} \rightarrow {h_i},
%  \sum_{j=1; j \ne i}^p \frac{\lambda_j\lambda_i}{(\lambda_i - \lambda_j)^2} .
  \label{eq:error}
\end{equation}
where
\begin{equation}
  h_i \stackrel{\Delta}{=} \sum_{j=1; j \ne i}^p \frac{\lambda_i \lambda_j}{(\lambda_i - \lambda_j)^2}
  \label{eq:hi}
\end{equation}
{and} $\lambda_1,\ldots,\lambda_p$ are the  population eigenvalues {of  $\bC$} (which we assume to be simple and in ascending order). 
%$\tilde{\lambda}_1,\ldots,\tilde{\lambda}_p$  and 
{One important application of the asymptotic result \eqref{eq:error} is that it provides an estimate for}
%These estimates can be used to bound 
the expected residual error between the sample and the population eigenvectors % $\bu_i - \tu_i$
% in expectation, asymptotically 
for large $n$, %using 
\begin{equation}
  \E{\|\bu_i - \tu_i \|^2} \approx \frac{1}{n}  {h}_i.
\label{eq:expectdelta}
\end{equation}
The usefulness of \eqref{eq:error}--\eqref{eq:expectdelta}, however, is limited by the fact that $h_i$ requires knowledge of all $p$ eigenvalues, which can be problematic---even computationally infeasible---when $p$ is large. 
%Particularly, solving for $p$ eigenvalues can become computationally infeasible for large matrices. 
Moreover, the values $\{\lambda_i\}$ are typically unknown for empirical data, and in practice one often approximates \eqref{eq:hi} using $\tilde{\lambda}_i\approx \lambda_i$, where $\{\tilde{\lambda}_i\}$ are the corresponding sample eigenvalues of $\tilde{\bC}$. 
}

{Thus motivated,  we study \eqref{eq:error}--\eqref{eq:expectdelta} for the limit of large $p$, seeking to avoid the computation of $p$ distinct eigenvalues by considering situations in which the right-hand side of \eqref{eq:hi} converges with increasing $p$. 
Defining such an extension, however, comes with several complications. One difficulty is that by allowing $p$ to increase, one ceases to study a single population covariance matrix $\bC$, but instead studies a sequence of population covariance matrices of growing size $p$.
%, each estimated by a corresponding sample covariance matrix $\tilde{\bC}$. 
One must therefore make an assumption about the origin of these population covariance matrices, and herein we assume they are drawn
%This requires one to assume that the population covariance matrices (in addition to the sample covariance matrices) are sampled 
from a matrix ensemble.
Moreover, because $\lim_{p\to\infty} \|\bu_i - \tu_i\|^2 $ is identical for any fixed $i$ (i.e., since $i/p\to0$ for  fixed $i$), we find it more interesting to study the $p\to\infty$ limiting behavior of \eqref{eq:error}--\eqref{eq:expectdelta} for fixed $\lambda_i=\lambda$, examining the associated eigenvectors $\{\tu_i\}$ across the ensemble. (Note that the index $i$ can vary from one population covariance matrix to another.)
%we consider $h_i$ for a neighborhood, $\{h_i: |\lambda_i-\lambda|<\delta\}$ for some small $\delta>0$. 
Finally, in this research we will assume  \eqref{eq:error} as a starting point---that is, we assume $n\to\infty$ much faster than $p\to\infty$. 
In Sec.~\ref{sec:discussion}, we study the necessary relative scaling behavior of $p$ and $n$,
finding that $n=\mathcal{O}(p^2)$ is a necessary relative scaling for the ensemble of population covariance matrices that we study.
% in Sec.~\ref{sec:discussion} and point out that it would be interesting in 
%Future research to extend \eqref{eq:error} by allowing $p$ to grow sufficiently slowly with $n$. 
}

{The first main contribution of this paper is an}
%The main contributions of this paper are 
asymptotic $p\to\infty$ estimate, {$\hat{h}_i \approx h_i$, for}
%of the quantities 
%\begin{equation}
%  h_i \stackrel{\Delta}{=} \sum_{j=1; j \ne i}^p \frac{\lambda_i \lambda_j}{(\lambda_i - \lambda_j)^2},
%  \label{eq:hi}
%\end{equation}
when the {population} eigenvalues {$\{\lambda_i\}$}
%$\lambda_1,\ldots,\lambda_p$ 
are distributed according to a known limiting $p\to\infty$ spectral density $\rho(\lambda)$. 
The idea of taking advantage of existing {\it a priori} knowledge about the spectral density $\rho(\lambda)$ has led to novel insights and improved inference for covariance analyses (\cite{Bickel2008,Lam2009,ledoit2011eigenvectors}). In practice, the probability distribution $\rho(\lambda)$ can be estimated from empirical data, or can sometimes be derived analytically (\cite{Mehta1991,Kuhn2008}). A situation of particular interest is when the covariance matrix follows a graphical (i.e., network-based) model in which complex network properties can give rise to different spectral densities (c.f., \cite{Chung2003,Goh2001,Zhang2004,Farkas2001,Dorogovtsev2003,peixoto2013eigenvalue,benaych2011eigenvalues,taylor2017super,taylor2016enhanced}).

{Note that values $h_i$ given by \eqref{eq:hi} depend on the consecutive 
%eigengaps $s^\pm  \stackrel{\Delta}{=} |\lambda_i-\lambda_{i\pm1}|$ for the population eigenvalues, 
%In order to estimate the size of the $h_i$ in \eqref{eq:hi}, we further use an approximation to the joint probability distribution $J(s^-,s^+)$ of the 
right and left  eigengaps around each eigenvalue $\lambda_i$,
\begin{equation}
  s_i^+ \stackrel{\Delta}{=} \lambda_{i+1} - \lambda_i \; \text{and}\;
  s_i^- \stackrel{\Delta}{=} \lambda_{i} - \lambda_{i-1}.
\end{equation}
In the context of quantum physics, these eigengaps are often referred to as level spacings since the eigenvalues typically represent energy levels (\cite{guhr1998random}).
Herein, we assume the population covariance matrices have eigengap statistics consistent with the Gaussian orthogonal ensemble (GOE) of random matrices, thereby allowing us to take advantage of existing theory for GOE eigengap statistics. In particular, we leverage the Wigner surmise (\cite{wigner1958distribution,wigner1993class})
\begin{equation}\label{eq:surmise}
    %P_s(s) \approx 
    P(s) = \frac{\pi p^2\rho^2(\lambda)}{2} s ~\text{exp}\left({-\frac{\pi p^2\rho^2(\lambda)}{4} s^2} \right)
\end{equation}
for $s\in\{s_i^{+},s_i^{-}\}$, which is a celebrated result  obtained by Eugene Wigner in the 1950's to describe the distribution of eigengaps for GOE matrices of size $p=2$. Equation \eqref{eq:surmise}  has had an enormous impact in physics (\cite{abul1999wigner,brody1973statistical,ellegaard1995spectral,pimpinelli2005evolution,schierenberg2012wigner,shklovskii1993statistics}) and economics (\cite{plerou2002random,akemann2010universal}). It was originally introduced as a `surmise' because it was believed to offer an accurate approximation to the eigengaps for large GOE matrices. Remarkably, over the last 5 decades there has been considerable numerical support validating the approximation's accuracy for large GOE matrices in which $p\gg2$. Moreover, \eqref{eq:surmise} has been observed to accurately predict the eigengap distribution for numerous empirical covariance matrices describing real-world datasets (\cite{plerou2002random,akemann2010universal}). 

Our second and third main contributions  utilize an extension to the Wigner surmise that approximates the joint distribution $J(s^-,s^+)$ %\approx P_{s^+,s^-}(s^+,s^-)$ 
and was derived for GOE matrices of size $3\times 3$ (\cite{Herman2007}). 
While developing theory based on such approximations  introduces error into our analysis, as we shall show, the simplicity of these \emph{surmises} allows us to make insights that may otherwise be unobtainable. 
%Finally, we note that because we are interested in computing estimates of $h_i$ in \eqref{eq:hi} when both $p\rightarrow \infty$ and $n\rightarrow \infty$, we expect that the gaps between any two eigenvalues, $s_i^+$ and $s_i^-$, will go to zero. 
%making such an approximation has two important consequences. 
For example,  it is easy to show from  \eqref{eq:surmise} that {an} expected eigengap 
%$s_i^{\pm}$ %between two eigenvalues $\lambda_i$ and $\lambda_{i+1}$ 
should have size 
\begin{equation}
   {\mathbb{E}[s^\pm] = \mathcal{O}\left( \frac{1}{p \rho(\lambda_i)} \right)}
   \label{eq:limit_s}
\end{equation}
as $p \rightarrow \infty$. 
See also (\cite{pastur11}, p.~16)  in which this scaling is obtained as the ``typical spacing unit'' for a random matrix with convergent spectral density.
%
%As a result, many of our results will be expressed in terms of this limiting eigengap.}
%
%{Second, the Wigner-surmise assumption inherently introduces an approximation error, and understanding this error remains an important open problem for the general  study of GOE matrices in general as well as a interesting direction for future work to extend this research; nevertheless, by making such an assumption, we are able to leverage the simplicity of these surmises to obtain useful predictions regarding the large $p$ behavior of $h_i$ and the necessary relative scaling between $p$ and $n$. 
%}
%
In this work, we use \eqref{eq:limit_s} to study the large-$p$ scaling behavior for $h_i$  as well as  the necessary relative scaling between $p$ and $n$.
Our approach involves introducing and estimating a probability density function $f_H(h)$ of $h_i$ in \eqref{eq:hi}, which describes the distribution of $h_i$ (keeping $\lambda_i$ fixed) across the population-covariance-matrix ensemble.
% from which the population covariance matrix is drawn (e.g., the GOE).
We obtain  estimates for $f_H(h)$ in terms of $\lambda_i$, $\rho(\lambda_i)$ and $p$, which
 %Such approximations 
are of great consequence, because they describe the expected uncertainty associated with  sample eigenvectors
%, 
%using a known limiting spectral density $\rho(\lambda)$, 
%thereby providing an estimate of eigenvector uncertainty 
across the  ensemble of population covariance matrices associated with $\rho(\lambda)$. That is, our second and third main results offer estimates for the expected eigenvector error $\E{\|\bu_i - \tu_i \|^2} $ that neither require a covariance matrix nor its eigenvalues. Importantly, the ensemble-based approach that we develop herein provides a new direction for uncertainty quantification of empirical eigenvectors that is  scalable for high-dimensional (large $p$) data.}
%---these  assumption regarding the spectral density of eigenvalues $\rho(\lambda)$.}

The paper is organized as follows. We state our main results in section \ref{sec:theory}.  In section \ref{sec:num}, we provide numerical simulations to {support these results}. In section \ref{sec:discussion}, we describe conditions in which \eqref{eq:error}, and thus our main results, are valid. The Appendix contains the derivations of our main results.

%________________________________________________________________
\section{Main results\label{sec:theory}}
%________________________________________________________________

In this section, we present asymptotic ($n\to\infty$ and $p\to\infty$) approximations for the expected residual error {  $\E{\|\bu_i - \tu_i \|^2}$} of sample eigenvectors as well as their distribution across an ensemble of population covariance matrices. 
We first provide preliminary discussion in section \ref{sec:assumptions}. 
{In section \ref{sec:main1}, we present main result 1, which provides a $p\to \infty$ estimate for the right-hand side of \eqref{eq:hi} using the assumption that the distribution of population eigenvalues weakly converges to a spectral density $\rho(\lambda)$.
In sections \ref{sec:main2} and \ref{sec:main3}, we present main results 2 and 3, which additionally assume the distribution of eigengaps for population covariance matrices is the same as that for the GOE random-matrix ensemble.}

%
%provides a $p\infty$ estimate for the right-hand side of 
%
%--\ref{sec:main3}, we present our three main results. 

\subsection{{Model specification} and assumptions\label{sec:assumptions}}
%\subsection{Definitions, notations and assumptions\label{sec:assumptions}}
%________________________________________________________________

{We consider a sequence of population covariance matrices (each denoted $\bC$) of growing size $p\to\infty$ such that each is drawn from a random-matrix ensemble. Let $\{\lambda_i\}_{i=1}^p$ and $\{\bu_{i}\}_{i=1}^p$ denote the population eigenvalues and corresponding eigenvectors, respectively, for each $\bC$. We 
%and $\rho_p(\lambda) = p^{-1}\sum_p \delta_{\lambda_i}(\lambda)$ denote the corresponding spectral density, 
 make the following two assumptions regarding the eigenspectra for the matrix ensemble.}

\begin{assume}\label{asssume1}
{We assume that the population eigenvalues $\{\lambda_i\}$ are simple and that the spectral density $\rho_p(\lambda) = p^{-1}\sum_{i=1}^p \delta_{\lambda_i}(\lambda)$ 
%are identically distributed with the limiting $p\to\infty$ 
weakly converges as $p\to\infty$ to a limiting spectral density $\rho_p(\lambda)\to \rho(\lambda)$
%that is , and we assume weak convergence {at a rate that is sublinear, i.e., $|\int f(\lambda)[ \rho_p(\lambda) - \rho(\lambda)] d\lambda | = \mathcal{O}(p^{-\gamma})$ with $\gamma>1$ for any bounded continuous function $f(\lambda)$. (Here $%\rho_p(\lambda) = p^{-1}\sum_p \delta_{\lambda_i}(\lambda)$.)} We further assume that $\rho(\lambda)$ 
that} has compact support $[\lambda_\text{min},\lambda_\text{max}] \subset \R^+$ and is continuous and differentiable on the interior of its support, $(\lambda_\text{min},\lambda_\text{max})$.
%
%
%t{We further assume a sublinear rate of weak convergence such that 
%\begin{equation}
%\int_{\lambda_\text{min}}^{\lambda_\text{max}} f(\lambda)[ \rho_p(\lambda) - \rho(\lambda)] d\lambda  = \mathcal{O}(p^{-\gamma})
%\end{equation}
% for any bounded continuous function $f(\lambda)$ as $p\to\infty$ and $\gamma>1$.}

\end{assume}

Many ensembles of symmetric random matrices satisfy Assumption~\ref{asssume1} (see \cite{Mehta1991,Kuhn2008,anderson03}) including, for example, those described by the semi-circle law (\cite{pastur11}, see sections 2 and 6). For some applications, it may also be beneficial to posit a parametric model for $\rho(\lambda)$, which can be estimated for small $p$ and $n$ and extended to the entire dataset. \\

\begin{assume}\label{asssume3}
We assume that the joint probability distribution $J(s^-,s^+)$ of the left and right gaps, {$s^\pm_i = |\lambda_i - \lambda_{i\pm 1}|$,}  around each eigenvalue $\lambda_i$ {is given by the following generalized Wigner surmise for the GOE:}
%can be accurately approximated by
\begin{equation}
 J(s^- ,s^+) \approx 
  \frac{3^7\left[p \rho(\lambda_{{i}})\right]^5}{32 \pi^3}
  \left[s^+s^- (s^+   +   s^- )\right] 
 \exp\left(-\frac{\left[ 3p\rho(\lambda{_{i}}) \right]^2 }{4\pi}
  \left[(s^+ )^2 +(s^-)^2+s^+  s^- \right]
  \right). 
\label{eq:joint} 
\end{equation}
\end{assume}

The joint distribution given by \eqref{eq:joint} was derived in \cite{Herman2007}  {[see equation (15)]} using $3 \times 3$ GOE matrices %from the Gaussian orthogonal ensembles (GOE) 
and can be {constructed} as a generalization of the Wigner surmise {[see \eqref{eq:surmise}]}, which approximates marginal distributions for $J(s^- ,s^+) $.
% and can be derived for $2 \times 2$ GOE matrices. 
{In addition to establishing the scaling $\mathbb{E}[s^\pm] = \mathcal{O}\left(  {1}/{p \rho(\lambda_i)} \right)$ [see \eqref{eq:limit_s}],} assumption \ref{asssume3} {also} implies that the $p$ population eigenvalues are {simple (i.e., distinct),} $\lambda_1 < \cdots < \lambda_p$, and is akin to the ``level repulsion of eigenvalues'' observed in large random matrices that states that the {eigengap} probability is 0 for $s^{\pm}=0$ (\cite{bourgade14}). 
%Moreover, this assumption establishes $\mathcal{O}(1/p\rho(\lambda))$ as the asymptotic scaling for the expectation of $s^{\pm} $.

As shown in \cite{Herman2007}, \eqref{eq:joint} gives a very good approximation to the exact distribution, {which may be expressed as an infinite dimensional integral and can be approximated using numerical integration and Toeplitz determinants. 
We note in passing, that it would be interesting in future research to connect this approach  to \eqref{eq:error}--\eqref{eq:expectdelta}; however, it is unclear whether or not this approach would allow for the type of results as we present here. In contrast, while the surmises [i.e., \eqref{eq:surmise} and \eqref{eq:joint}] introduce an approximation error---which  remains an open topic of great interest in random matrix theory---their simplicity allows us to gain insights that may be otherwise unobtainable.
%be difficult or impossible to obtain.
% one could take this approach to study the limiting $p\to\infty$ behavior of $h_i$, we do not find the existing theory to be amenable for obtaining  closed-form approximations to $h_i$ allowing important insights such as scaling behaviors. In contrast,
%
In addition, as demonstrated in our numerical simulations {(see Figures \ref{rho} and \ref{joint_fig})}, \eqref{eq:joint} provides a good approximation for the ensemble of {population} covariance matrices that we study.}
% (see Figure~\ref{joint_fig}). 

{
For each population covariance matrix $\bC$,
we consider a sample covariance matrix
%We consider the $p \times p$ sample covariance matrix 
%{
$\tilde{\bC}= n^{-1} (\bX -\mathbb{E}[\bX ])(\bX-\mathbb{E}[\bX])^T$ 
 computed from $n$ observations, $\bx_1,\ldots, \bx_n$, of a random vector $\bx \in \R^p$ {and $\bX = [\bx_1,\dots,\bx_n]$}. We denote by $\tlamb_1 \leq \cdots \leq  \tlamb_p$ the $p$ sample eigenvalues of $\tilde{\bC}$, and $\tu_1,\cdots,\tu_p$ the corresponding sample eigenvectors. 
We assume each sample covariance matrix $\tilde{\bC}$ is Wishart-distributed around ${\bC}$, as is the case when $\bx$ follows a multivariate-Gaussian distribution (\cite{anderson03}).
 }
 %Likewise, $\lambda_1\leq \cdots \leq \lambda_p$ and $\bu_1,\cdots,\bu_p$ denote the corresponding population eigenvalues and associated eigenvectors, respectively, of ${\bC}$.}
%

%{Our main results will require the following technical assumptions.}\\

%________________________________________________________________
\subsection{Main result  1: Estimate of $h_i$ for large $p$}\label{sec:main1}
%________________________________________________________________ 

{We may now present our first main result, an estimate  {$\hat{h}_i $ for $h_i$ [see \eqref{eq:hi}]} for high-dimensional (large $p$) covariance matrices.}
{Under Assumption \ref{asssume1}, } 
%the limit of large $p$ and large $n$ , 
we find the following {asymptotic} $p\to\infty$ approximation for $h_i$,
\begin{equation}
  \hat{h}_i  = \lambda_i^2 \left[
        \left(
        \frac{1}{(s_i^-)^2}  + \frac{1}{(s_i^+)^2}
    \right)
    + p \rho(\lambda_i)  
    \left(
        \frac{1}{s_i^-}  + \frac{1}{s_i^+}
    \right)
\right]. \label{eq:main_approx} 
\end{equation}
See Appendix~\ref{sec:proof_main1} for the derivation.

We can explain the role of the different terms in \eqref{eq:main_approx} as follows. The left term in the squared brackets approximates the terms in \eqref{eq:hi} that involve the nearest neighbor eigenvalues of $\lambda_i$, which are respectively located at $\lambda_{i-1}=\lambda_i - s_i^{-}$ and $\lambda_{i+1}=\lambda_i + s_i^{+}$ and dominate the summation in \eqref{eq:hi} when $p$ is large. This term does not require the knowledge of the probability distribution $\rho(\lambda)$. The second term accounts for the remaining terms in \eqref{eq:hi}, which involve the remaining eigenvalues, $\{\lambda_j : |j-i|>1\}$. Finally, as explained {for \eqref{eq:limit_s}, since}
%, we recall that 
$1/ p \rho(\lambda_i)$ {is} the same order as the {expected} gap between two {population} eigenvalues (\cite{pastur11}, {p.~16}), %and %therefore 
all terms in the right-hand side of \eqref{eq:main_approx} can potentially obtain similar magnitudes.

{Note that main result 1 does not depend on $n$. It  approximates $h_i$ for population covariance matrices drawn from a matrix ensemble with a convergent spectral density $\rho(\lambda)$. By combining \eqref{eq:main_approx}  with \eqref{eq:expectdelta}, we obtain a large-$p$ approximation to the expected error of a sample eigenvector $\tu_i $:}
%We note that \eqref{eq:main_approx} allows one to approximate 
\begin{equation}
  \E{{\|\bu_i - \tu_i \|^2}}\approx \hat{h}_i/n\label{eq:PUG}.
\end{equation}
%for the sample eigenvector $\tu_i $ associated with eigenvalue $\lambda_i$ of a large covariance matrix (i.e., $p\gg1$) using 
{Importantly, this result uses} knowledge of the population eigenvalue distribution $\rho(\lambda)$ and {nearest-neighbor eigengaps} $s_i^\pm$. Thus, it does not require knowledge of the full set of sample eigenvalues $\{\tlamb_i\}$, {and is therefore  scalable for high-dimensional (large $p$) data.} However, we stress that approximation {\eqref{eq:PUG}} is valid only when $p$ and $n$ are sufficiently large, which we will explore in section \ref{sec:discussion}.

\clearpage

%________________________________________________________________
\subsection{Main result  2: Estimate of the probability density of {$h_i$ across matrix ensemble}\label{sec:main2}}
%________________________________________________________________ 

{Equation \eqref{eq:main_approx} gives an asymptotic estimate for $h_i$ that uses $\lambda_i$, $s_i^+$ and $s_i^-$, which can be calculated for a given population covariance matrix $\bC$ (or estimated from $\tilde{\bC}$). We now turn our attention to studying the distribution of $h_i$---denoted  $f_H(h_i)$---across all the population covariance matrices in the matrix ensemble for which $\lambda_i=\lambda$ is an eigenvalue.  That is, we consider fixed $\lambda_i$ and approximate $f_H(h_i)$ by allowing $(s_i^+,s_i^-)$ to be distributed according to assumption \ref{asssume3}. Note that, once $\lambda$ is fixed, the index $i$  can differ from one population covariance matrix to another, and so from here on we will drop the subscript $i$ whenever describing the distribution of a variable across the matrix ensemble.}

%\todo{carefully define distribution of $\hat{h}_i$. }

We again consider the case where $p$ is large, and we study the limit of the distribution of {$\hat{h}$ given by \eqref{eq:main_approx}} for $p \rightarrow \infty$. By combining \eqref{eq:main_approx} {with} 
%the expression of the joint probability for the eigenvalue gap given by 
{\eqref{eq:joint}}, we obtain the following semi-analytical expression for the limiting probability density of the approximation $\hat{h}$ to $h$, 
\begin{equation}
   f_H(h) =  - \int_{s^0(h)}^{\infty} J\left(s^*(h,s^+),s^+\right) \frac{\partial  s^*(h,s^+)}{\partial h}ds^+,
   \label{eq:main_density}
\end{equation}
where the variables $s^0(h)$ and $s^*(h,s^+)$ depend on the eigenvalue $\lambda$ around which $h$ is computed, and are respectively given by 
\begin{equation} 
s^0(h)  = \frac{\lambda^2 p \rho(\lambda)}{2h}
\left \{ 
1 + \sqrt{1 + \frac{4h}{\left[\lambda p \rho(\lambda)\right]^2}} 
\right\},
\end{equation}
and 
\begin{equation}
  s^*(h,s^+) = \lambda^2 p \rho(\lambda)  
  \frac{
    1 +  \sqrt{ 1 + \frac{\displaystyle 4}{\left[\displaystyle \lambda p \rho(\lambda)\right]^2}
       \left(
         h - \frac{\displaystyle \lambda^2}{\displaystyle (s^+)^2} - \frac{\displaystyle\lambda p \rho(\lambda) }{\displaystyle s^+} 
       \right)
    }
  }
  {
    2\left(h - \frac{\displaystyle \lambda^2}{\displaystyle (s^+)^2} -
      \frac{\displaystyle \lambda p \rho(\lambda) }{\displaystyle s^+}
    \right)
  }.
\end{equation}
See Appendix~\ref{sec:estimate2} for the derivation.

The significance of \eqref{eq:main_density} stems from the fact that it allows one to approximate the distribution of $h$ and therefore the distribution of expected eigenvector errors using the approximation \eqref{eq:PUG}, which again assumes sufficiently large $p$ and $n$. Specifically, $f_H(h)$ estimates the distribution of expected residual error across the covariance-matrix ensemble associated---that is, as opposed to \eqref{eq:error}, which is an estimate for a single covariance matrix from the ensemble.

%________________________________________________________________
\subsection{Main result  3: Asymptotic behavior of $f_H(h)$ for large $h$}\label{sec:main3}
%________________________________________________________________ 
Keeping $\lambda$ and $\rho(\lambda)$  fixed, in the limit when the left gap, $s^-$, or right gap, $s^+$, goes to zero, then $h$ goes to
infinity, and we find the following scaling behavior for the probability density function 
\begin{equation}
   f_H(h)  = \mathcal{O}\left(\frac{p^2}{h^2}\right). 
   \label{eq:density_largep}
\end{equation}
See Appendix~\ref{sec:limit} for the derivation. 
%Therefore, $f_H(h)$ scales like $h^{-2}$ for large $h$. 

{We point out that the limit of large $h$ is especially interesting because it corresponds to the case where the sample estimates $\tilde{\bm{u}}_i$ of the eigenvectors $\bm{u}_i$ are the least accurate. The observation that $f_H(h)$ has a power-law decay for large $h$ implies that the error associated with sample eigenvectors is very heterogeneous, and one should expect situations in which the error $\E{{\|\bu_i - \tu_i \|^2}}\approx  {h}_i/n$ is small for many $\tu_i$, but on the rare occasions in which $s_{i}^{+}$ and/or $s_i^-$ are  small, $\E{{\|\bu_i - \tu_i \|^2}}$ can be orders-of-magnitude larger. 
Later, in Sec.~\ref{sec:discussion} we show that  $n\ge h_i/2$ is a necessary (but not sufficient condition) for \eqref{eq:expectdelta} to offer an accurate approximation. Because $h_i=\mathcal{O}(p^2)$, we find $n=\mathcal{O}(p^2)$ to be a necessary (but not sufficient) relative scaling for \eqref{eq:expectdelta} when $n,p\to\infty$. This should be further explored in future work.
}

%\clearpage

%___________________________________________________________
\section{Numerical {validation of main results}}\label{sec:num}
%___________________________________________________________
We now report the results of numerical experiments to validate the theoretical predictions given by the main results described in section \ref{sec:theory}. 
{In Sec.~\ref{sec:num_Laplacian}, we}
describe the {matrix} ensemble {for population} covariance matrices used for these experiments.
{In Sec.~\ref{sec:num_2}, we support main result 1.
We support main results 2 and 3 in Secs.~\ref{sec:num_4} and \ref{sec:num_6}, respectively.
}

%___________________________________________________________
\subsection{{Population covariance matrix ensemble: Laplacians of $k$-regular graphs}\label{sec:num_Laplacian}}
%Setting:  Laplacian of $k$-regular graphs\label{sec:num_Laplacian}}
%___________________________________________________________

{
We seek to study the error of sample eigenvectors in the limit of large $p$ and $n>p$, focusing on the scenario in which the population covariance matrices are drawn from a matrix ensemble satisfying assumptions \ref{asssume1} and \ref{asssume3}. All covariance matrices must also be positive semi-definite (\cite{anderson03}) so that $\lambda_i\ge0$ for all $i$. In addition, to help mitigate the computational cost of studying the eigenspectra for high-dimensional (large $p$) covariance matrices, we would like to study a sparse random-matrix ensemble in which most matrix entries are zero. 
%For example, the GOE random-matrix ensemble is often numerically studied through the stochastic generation of spectrally equivalent random tridiagonal matrices (\cite{Mehta1991}). We pursue a different direction by allowing the pattern of matrix sparsity to also be random; 

Thus motivated, we study a graphical model (\cite{Hastie2009}).}
{We let the population covariance matrices be given by the unnormalized---also called combinatorial (\cite{Bapat10})---Laplacian matrices\footnote{{Any Laplacian matrix $\bC$ is positive semi-definite: ${\bf v}^T \bC {\bf v} \ge0$ for any vector ${\bf v}$.
%\todo{something here first???}
%We consider the following graphical model for {the population} covariance matrices {in which they} 
%covariance matrices $\bC$ 
Moreover, Laplacian matrices arise for many types of random processes on graphs and are related, for example, to the autocovariance matrices of random walks on graphs (\cite{Delvenne2010}).
We also note that a Laplacian matrix $\bC$ can be written as
%take the form 
%\begin{equation}
$  \bC = \bX \bX^T,$
%  \label{laplacian}
%\end{equation}
where $\bX$ is a random incidence matrix that describes the connectivity of a random graph $G=(V,E)$, with an arbitrary orientation of the edges. {Each} entry $x_{e,v}$ of $\bX$  {can take one of three values:} 
1 if $v$ is the head of the oriented edge $e$;
-1 if $v$ is the tail of the oriented edge $e$;
{or} 0 otherwise.
%given by
%\begin{equation}
%$x_{e,v} = 
%\begin{cases}
%  \mspace{14mu} 1 & \text{if $v$ is the head of the oriented edge $e$,}\\
%  -1 & \text{if $v$ is the tail of the oriented edge $e$,}\\
%  \mspace{14mu}  0 & \text{otherwise.}
%\end{cases}
%\end{equation}
%We generate the random incidence matrix for a $k$-regular graph using the configuration model (\cite{Newman2003}). 
%We note that the covariance matrix in (\ref{laplacian}) is equal to the (unnormalized) combinatorial Laplacian of a graph (\cite{Bapat10}). 
}} of random $k$-regular graphs, which we generate using the configuration model (\cite{Newman2003}).
For $k$-regular graphs with fixed $k\gg1$, the spectral density $\rho_p(\lambda)=\sum_i \delta_{\lambda_i}(\lambda)$ weakly converges as $p\to\infty$ to a semicircle distribution 
\begin{equation}\label{eq:mckays}
{\rho_p(\lambda) \to }\rho(\lambda)=
\begin{cases}
\frac{\displaystyle k \sqrt{4 (k -1)  - \lambda^2}}{\displaystyle  2 \pi (k^2 - \lambda^2)} ,& \text{if $|\lambda| \leq 2
  \sqrt{k-1},$}\\
0 ,& \text{otherwise.}
\end{cases}
\end{equation}
Equation \eqref{eq:mckays} is known as McKay's law (\cite{McKay1981}).
While McKay obtained  \eqref{eq:mckays} for fixed $k\gg1$, it also describes the case for increasing $k$, provided that $k$ grows sufficiently slowly with $p$
%$ accurate when 
%McKay's law was originally obtained as the limiting spectral distribution for fixed $k$; however,  
%One can also allow $k$ to {slowly} increase with $p$ 
(\cite{dumitriu2012sparse}). 

% and averaging their spectral densities.}

Numerous empirical Laplacian matrices 
%describing real-world datasets  
have been observed to give rise to eigengap statistics consistent with the Wigner surmise given by \eqref{eq:surmise} (\cite{plerou2002random,akemann2010universal}), and we therefore believe the extended surmise given by \eqref{eq:joint} will also be widely applicable.
Importantly, our assumption that $k\gg1$ ensures all graphs are strongly connected, which has been observed to be an important requirement for the eigengap statistics to behave similarly to that for the GOE (\cite{murphy2017anderson}).  Understanding the relation between eigengap statistics and graph topology remains an important open topic  (\cite{murphy2017anderson,taylor2017super}). 
In future work, it would be interesting to allow for graphs with more complicated structure---often called complex networks---and there is a large body of work exploring spectral densities for these graphs
(\cite{Chung2003,Goh2001,Zhang2004,Farkas2001,Dorogovtsev2003,peixoto2013eigenvalue,benaych2011eigenvalues,taylor2017super,taylor2016enhanced}).
}

\begin{figure}[t] 
\centering
\includegraphics[width=0.39\textwidth]{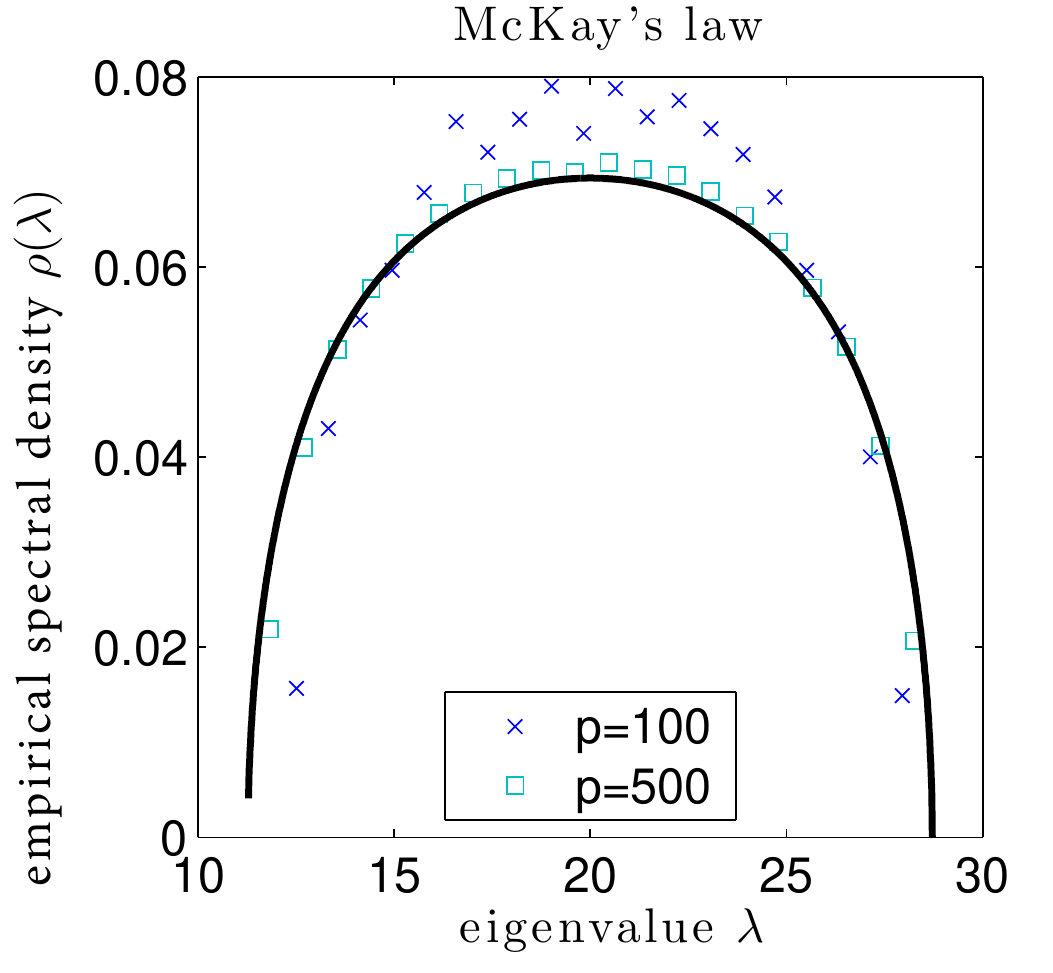}
\includegraphics[width=0.39\textwidth]{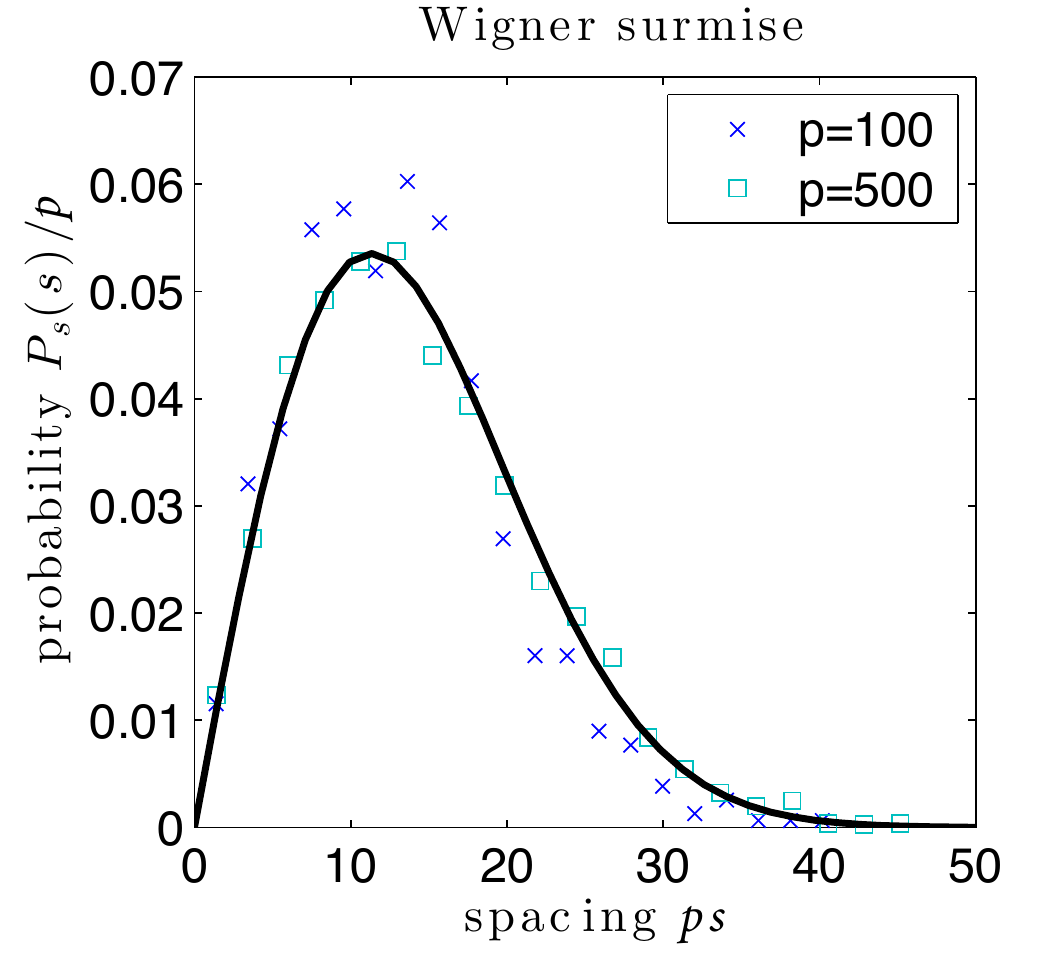}
\caption{
Left: Empirical spectral density {$\rho_p(\lambda)$ for {a} $k$-regular graphical model} converges toward McKay's law (\cite{McKay1981}) {given by \eqref{eq:mckays}} (black curve) in the limit  $p\to\infty$. 
Right: Distribution of normalized {eigengaps} $ \{p s_i^+\}$ for eigenvalues $|\lambda_i - 20 | < 1$ is well-described by the Wigner surmise {for GOE matrices, which is given by \eqref{eq:surmise}} (black curve). 
%(\cite{Gonzalez2008}). 
}
\label{rho}
\end{figure}
\begin{figure}[t] 
\centering
  \includegraphics[width=0.85\textwidth]{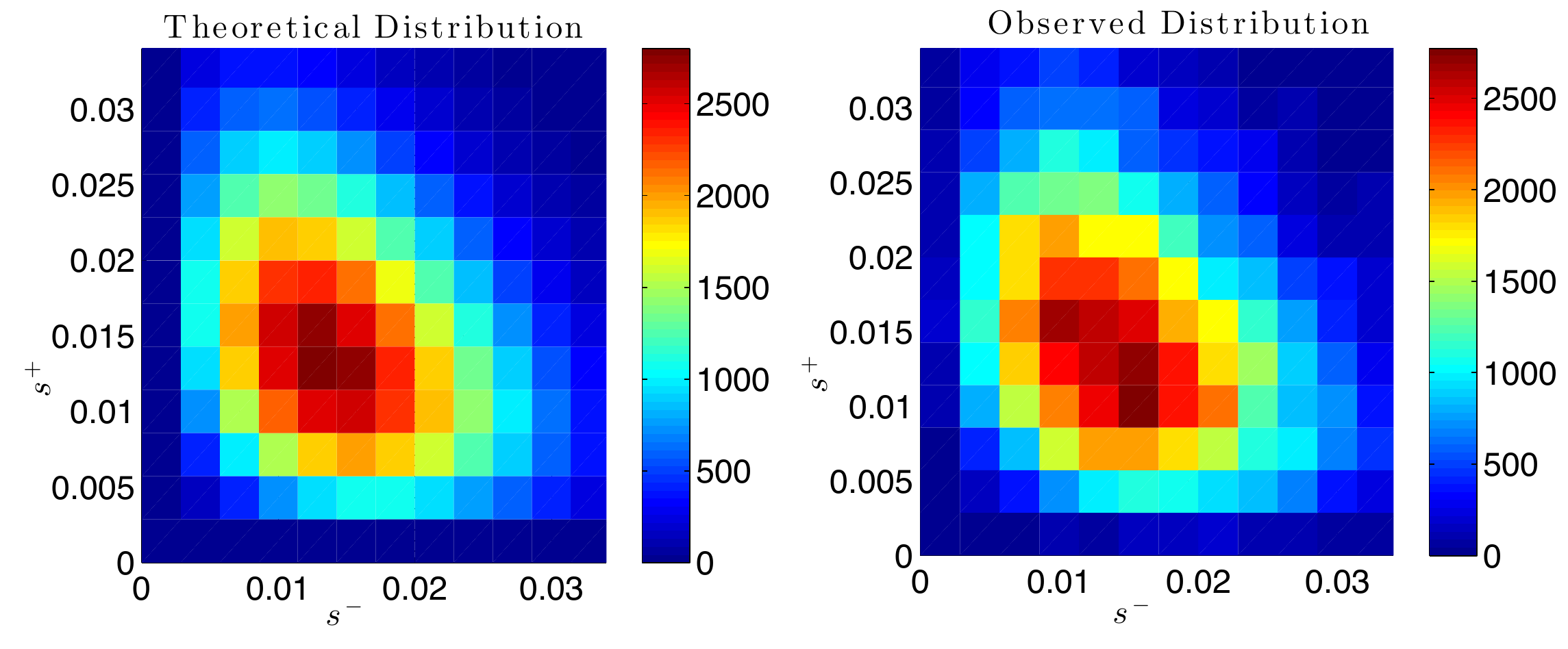}
\caption{Joint distribution of consecutive {eigengaps} $J(s^-,s^+)$: theoretical distribution (left) given by \eqref{eq:joint}, and numerically observed distribution (right). The color indicates the unnormalized counting measure. 
}
\label{joint_fig}
\end{figure}

{We now illustrate that this graphical model satisfies Assumptions 2.1 and 2.2.} Figure \ref{rho}-left displays the empirical spectral density computed over $50$  {population} covariance matrices {sampled from the graphical model}, using matrix sizes $p=100$ (blue crosses) and $p=500$ (cyan squares).  Note that $p$ indicates both the matrix size and the number of vertices in the {$k$-regular} graph. As $p$ increases from 100 to 500, we observe the convergence of the empirical spectral density toward {\eqref{eq:mckays}} (black curve).
%the limiting density given by McKay's law (\cite{McKay1981}) (black curve)
%
%\begin{equation}\label{eq:mckays}
%\rho(\lambda)=
%\begin{cases}
%\frac{\displaystyle k \sqrt{4 (k -1)  - \lambda^2}}{\displaystyle  2 \pi (k^2 - \lambda^2)} ,& \text{if $|\lambda| \leq 2
%  \sqrt{k-1},$}\\
%0 ,& \text{otherwise.}
%\end{cases}
%\end{equation}
%McKay's law was originally obtained as the limiting spectral distribution for fixed $k$; however,  one can also allow $k$ to {slowly} increase with $p$ (\cite{dumitriu2012sparse}). We also point out that in the numerical experiments to follow, we will estimate the distribution of eigenvalues directly from the data, rather than use McKay's law, as this approach would be more relevant for empirical data.
%
Figure \ref{rho}-right displays the empirical probability density of the normalized spacing $p s^+$ for the set of eigenvalues $\{\lambda_i\}$ such that $|\lambda_i - 20| <1$.  We used approximately $2 p \rho(\lambda)$ {eigengaps} to {estimate} the empirical densities. As expected, the {eigengap} distributions appear to be consistent with the Wigner surmise {given by \eqref{eq:surmise}} %(\cite{Gonzalez2008}) 
%\begin{equation}
%     P(s) \approx \frac{\pi p^2\rho^2(\lambda)}{2} s ~\text{exp}\left({-\frac{\pi p^2\rho^2(\lambda)}{4} s^2} \right),
%\end{equation}
(black curve). 
Note that the agreement improves with increasing  $p$. 
%\todo{move to earlier}
%
To gain further insight into the gap distribution, and to validate Assumption~\ref{asssume3}, we compared the unnormalized counting measure with the joint {eigengap} distribution $J(s^-,s^+)$  for the eigenvalues $\{\lambda_i\}$ such that $|\lambda_i - 20| <1$. Figure \ref{joint_fig}-left displays the level sets of $p^2 J(s^-,s^+)$ according to
\eqref{eq:joint} with $p=1,000$. Figure \ref{joint_fig}-right shows the unnormalized counting measure computed across $100$ covariance matrices of size $p=1,000$. These are in very good agreement.

{Before continuing, we need to make two clarifying points.
First, while we could use \eqref{eq:mckays} to test our approximations (see Sec.~\ref{sec:theory}), in the numerical experiments to follow,  we  instead estimate the limiting distribution of eigenvalues using the data, as this approach would be more relevant for empirical data. That is, we estimate $\rho(\lambda)$ by numerically computing the average spectral density of random population covariance matrices drawn from the graphical model for a given $p$.

Second, main results 2 and 3 describe the distribution $f_H(h_i)$ across the random-matrix ensemble from which population-covariance matrices are drawn. That is, we consider the distribution of $h_i$ associated with a particular eigenvalue $\lambda_i=\lambda$. However, if one fixes $\lambda$, then  one is confronted with an undersampling issue since it is  unlikely that the Laplacian of a randomly generated $k$-regular graph will have $\lambda$ as a particular eigenvalue.
%
%it is with very low probability that a random population covariance matrix will have $\lambda$ as an eigenvalue. Because generating large random matrices is computationally expensive, numerically studying $f_H(h)$ for fixed $\lambda$ requires that we address this undersampling issue. 
To overcome this issue,  we fix $\lambda$ and numerically study the distribution $f_H(h)$ for values $\{h_i\}$ associated with eigenvalues $ \{\lambda_i : |\lambda_i - \lambda|<\delta \}$ for small $\delta>0$. For each $p$, we choose $\delta$ to be sufficiently small so that the resulting distribution appears to not depend on $\delta$. We note that this represents a compromise between undersampling the random-matrix ensemble and the error introduced by allowing $\lambda_i$ lie within a small neighborhood (rather than remain fixed at $\lambda$).}

%___________________________________________________________
\subsection{Experimental validation of main result 1}\label{sec:num_2}
\begin{figure}[t] 
\centering
\includegraphics[width=0.44\textwidth]{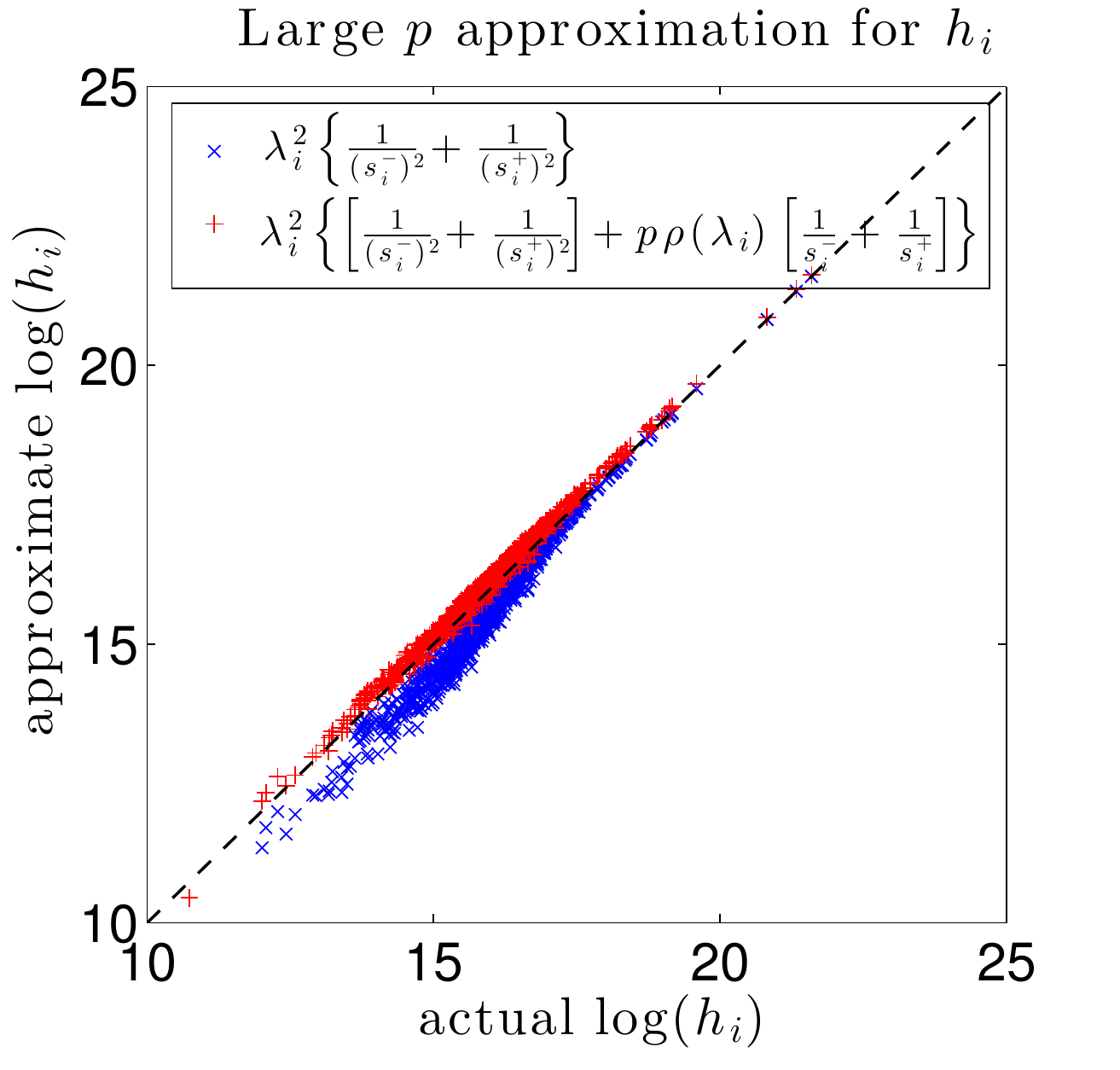}
\includegraphics[width=0.45\textwidth]{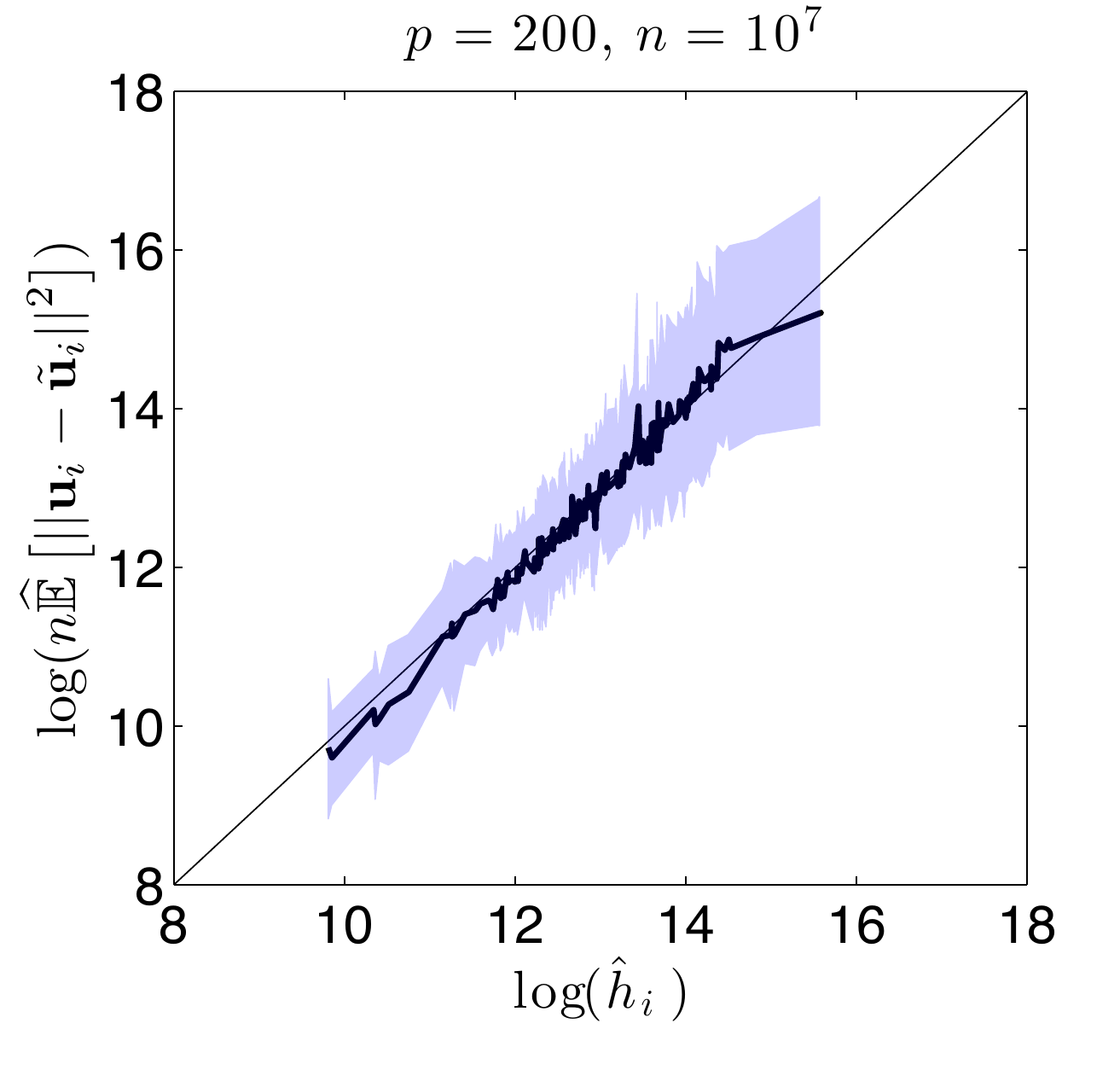}
\caption{{Support for main result 1.}
Left: Approximation $\widehat{h}_i$  
%of $h_i=\sum_{j\not=i}\lambda_i\lambda_j / (\lambda_i-\lambda_j)^2$ 
given by
\eqref{eq:main_approx} as function of the true $h_i$ {given by \eqref{eq:hi}. Results indicate $\hat{h}_i$ and  ${h}_i$ for a single population covariance matrix $\bC$ of size $p=100$ and $i\in\{1,\dots,p\}$.} We show $\widehat{h}_i$ with (red plus symbols) and without (blue crosses) the correction terms.
Right: Bootstrap estimate of the sample mean error, $n\wE{\|\bu_i - \tu_i \|^2}$, which is computed from $100$ samples from Wishart distribution $W(\bC,n)$ with $n=10^7$, versus approximation $\widehat{h}_i$ given by \eqref{eq:main_approx}. The mean is plotted by the black curve, and the standard deviation is shown in blue. See text for details.
}
\label{fig:accuracy2}
\end{figure}

%___________________________________________________________
We first compared the estimate $\widehat{h}_i$, given by \eqref{eq:main_approx}, with the true values of $h_i$, defined by \eqref{eq:hi}, for covariance-matrices ensemble described in section \ref{sec:num_Laplacian}. We considered graphs of fixed degree $k=20$ and $p=100$ vertices.
In Figure~\ref{fig:accuracy2}-left, we compare \eqref{eq:main_approx} with the true value of $h_i$ computed directly from the eigenvalues. The points lie close to the diagonal (dashed line), which validates the accuracy of the approximations. To illustrate the effect of the terms $p \rho(\lambda_i)\lambda_i^2/s_i^\pm$ in \eqref{eq:main_approx}, we plot our approximation with (red plus symbols) and without (blue crosses) these corrections. One can observe that these terms improve the estimate for small $h_i$ and have little effect for large $h_i$. This is expected since large 
$h_i$ corresponds to very small $s_i^\pm$. In this limit, the correction terms become negligible as $(s_i^\pm)^{-2}\gg (s_i^\pm)^{-1}$.

In the next experiment, we compare $\widehat{h}_i$ given by \eqref{eq:main_approx} with a bootstrap estimate of the mean sample error, $n\wE{\|\bu_i - \tu_i \|^2}$, for Wishart distribution $W(\bC,n)$. Specifically, we generated a population covariance matrix $\bC$ with $k=20$ and $p=200$. We then generated 100 random realizations $\widetilde{\bC}$ from  $W(\bC,n)$ with $n=10^7$.  Let $\{\bu_i\}_{i=1}^p$ be the eigenvectors of $\bC$.  For each random realization $\widetilde{\bC}$, we calculated its eigenvectors $\{\tu_i\}_{i=1}^p$ and computed the residual error $\bu_i - \tu_i$ between the sample eigenvectors and the population eigenvectors. We then computed a bootstrap estimate, $\wE{\|\bu_i - \tu_i \|^2}$, indicating the observed mean eigenvector error across the 100 realizations of $\widetilde{\bC}$.
In Figure~\ref{fig:accuracy2}-right, we plot the observed values $n\wE{\|\bu_i - \tu_i \|^2}$ versus our prediction given by \eqref{eq:PUG}. The mean is plotted in black, and the standard deviation is shown in blue. We note that the solid curves lie very close to the diagonal indicating the accuracy of \eqref{eq:PUG}. 

In these experiments, the sample size $n$ was chosen to be sufficiently large so that \eqref{eq:expectdelta} and \eqref{eq:PUG} are accurate. 
{Recall that \eqref{eq:error} is an asymptotic $n\to\infty$ limit for $  \E{n \|\bu_i - \tu_i\|^2} $, and we numerically observe that $n$ must be very large for the asymptotic result to provide an accurate approximation. }
We discuss in section \ref{sec:discussion} a simple and practical bound that can be used to choose appropriate values of $n$.

\begin{figure}[h] 
\centering
\includegraphics[width=0.45\textwidth]{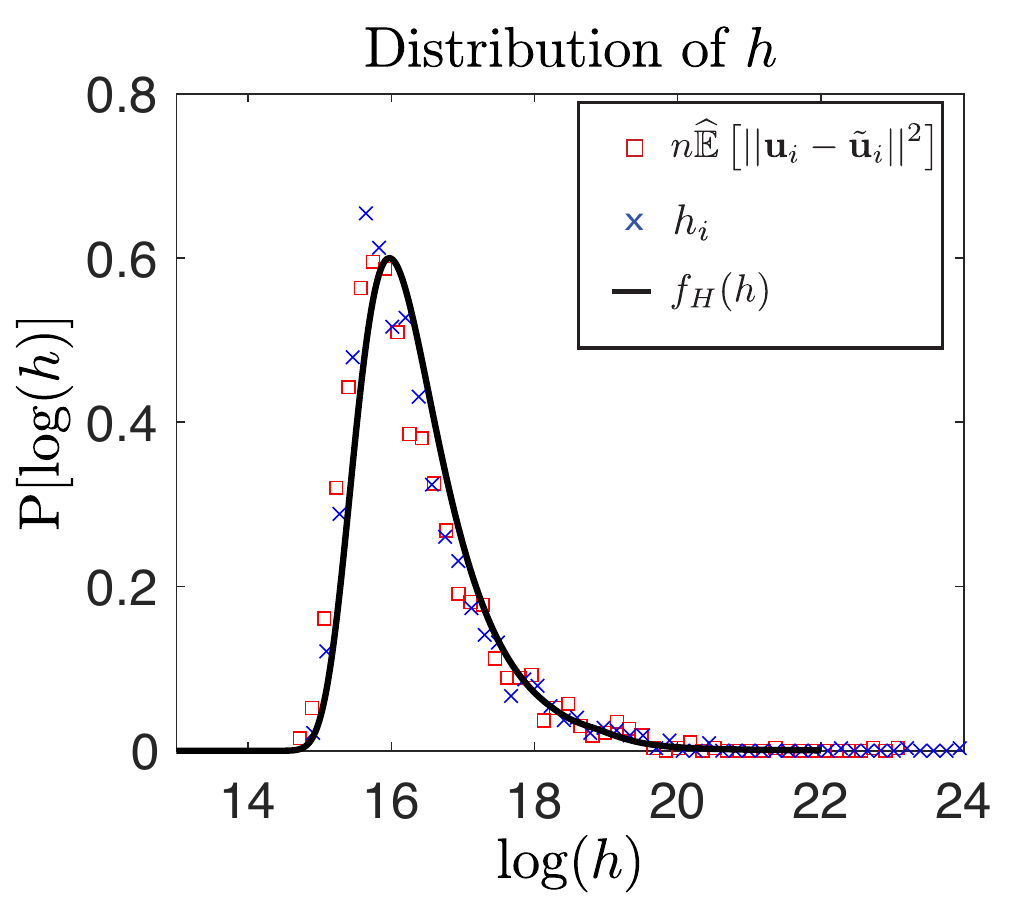}
\includegraphics[width=0.45\textwidth]{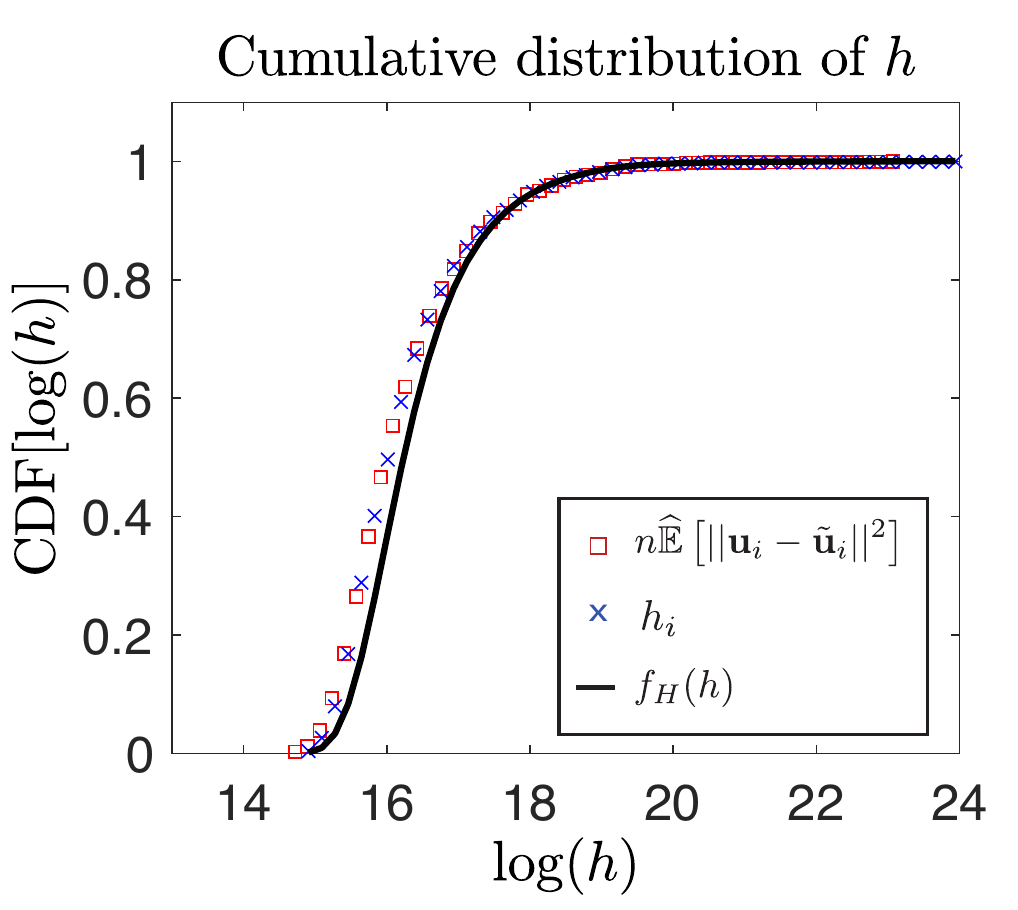}
\vspace{-.1cm}
\caption{{Support for main result 2.}
Accuracy of probability density function {$f_H(h)$ of $h$ given by \eqref{eq:main_density}} (left) and {its associated} cumulative density (right).
% of $h$ for $n=10^{10}$. 
We depict the following:
(black curves) semi-analytical expression of the probability distribution $f_H(h)$ given by \eqref{eq:main_density}; %and corresponding cumulative distribution;
(blue crosses, $\blue \times$): {an} empirically observed distribution of ${h}_i$ {computed for $\{h_i:|\lambda_i-\lambda|<1\}$ with $\lambda=20$}; 
(red squares, $\red \Box$): empirically observed distribution of bootstrap estimate, $n \wE{\|\bu_i - \tu_i\|^2}$.
}
\label{fig:accuracy4}
\end{figure}

%___________________________________________________________
\subsection{Experimental validation of main result  2 \label{sec:num_4}}
%___________________________________________________________

We now describe experiments that validate the second main result presented in section \ref{sec:main2}.  We confirm that the approximation of $f_H(h)$ given by \eqref{eq:main_density} is in good agreement with the empirical distribution of $h_i$. Furthermore, we show experimentally that $f_H(h)$ in \eqref{eq:main_density} also approximates the distribution of the expected residual error $\E{n{\|\bu_i - \tu_i \|^2}}$, provided that $n$ and $p$ are sufficiently large.

We generated 50 unweighted graphs of fixed degree $k=20$ and fixed size $p=1,000$. For each graph, we constructed the population covariance matrix, $\bC$, as explained in section \ref{sec:num_Laplacian}. For each $\bC$, we generated 10 sample covariance matrices $\widetilde{\bC}$ from Wishart distribution $W(\bC,n)$ with $n=10^{10}$. For each $\widetilde{\bC}$, we calculated its eigenvectors $\{\tu_i\}_{i=1}^p$ and computed the residual error, $\bu_i - \tu_i$, between the sample  and  population eigenvectors. {We consider all eigenvectors such that their associated eigenvalues satisfy $|\lambda_i-\lambda|<1$.} We then computed a bootstrap estimate, $n\wE{\|\bu_i - \tu_i \|^2}$, of the mean sample error for each $\bC$ using the 10 realizations of $\widetilde{\bC}$.

In Figure \ref{fig:accuracy4}-left, we {use} a  solid black curve {to represent} the semi-analytical expression of the probability distribution $f_H(h)$ given by \eqref{eq:main_density}. We plot its corresponding cumulative distribution in Figure \ref{fig:accuracy4}-right. We plot with blue crosses in both panels a numerically observed distribution of ${h}_i$, which we estimate using 50 covariances $\bC$ drawn from the graphical model described in section \ref{sec:num_Laplacian}. We plot with red squares an empirical distribution of bootstrap estimates, $n \wE{\|\bu_i - \tu_i\|^2}$. As expected, the probability density function $f_H(h)$ provides a good approximation of the empirical distribution of ${h}_i$ as well as the distribution of $n\wE{\|\bu_i - \tu_i \|^2}$ (that is, provided  $n$ and $p$ are both sufficiently large). However, we note that the distribution $f_H(h)$ is shifted slightly to the right.
{This is in agreement with Figure \ref{fig:accuracy4}, where one can observe that $\hat{h}_i$ typically  overestimates $h_i$ a very small amount (i.e., the red $+$ symbols tend to be just above the diagonal).}

%\clearpage 

\begin{figure}[h] 
\centering
\includegraphics[width=0.45\textwidth]{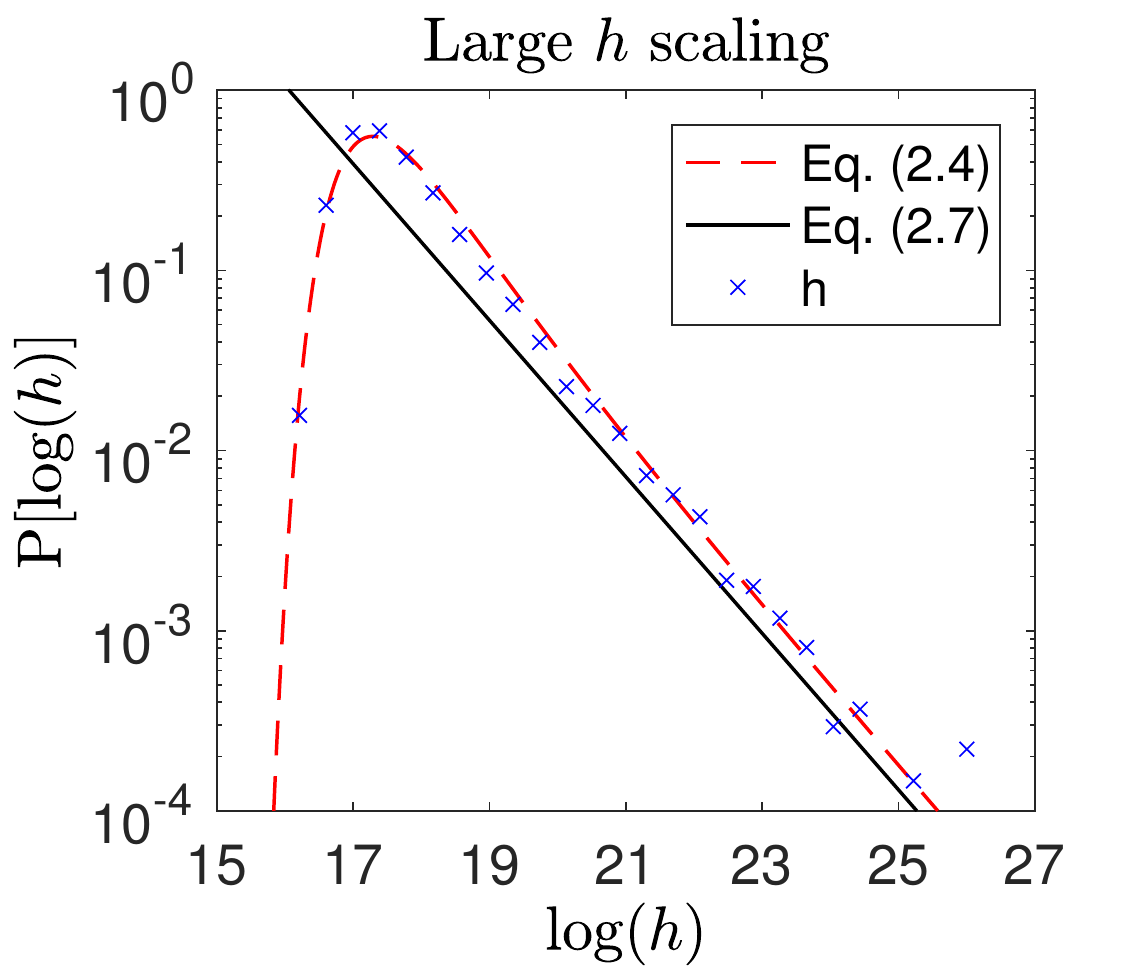}
\caption{{Support for main result 3.}
%Large-$h$ scaling 
{$\mathcal{O}(p^2/h^2)$ scaling of $f_H(h)$ in the limit of large $h$.} We plot P$[\log(h)]$ as a function of $\log(h)$: semi-analytical expression $f_H(h)$ computed from \eqref{eq:main_density} in red ({\red -~-}); limiting approximation to $f_H(h)$ for large $h$, given by \eqref{eq:density_largep} in black (--); empirical distribution, shown as blue crosses ($\blue \times$). 
}
\label{fig:accuracy6}
\end{figure}

%___________________________________________________________
\subsection{Experimental validation of main result  3} \label{sec:num_6}
%___________________________________________________________
We conclude with numerical validation of main result 3, $f_H(h)\varpropto h^{-2}$ for large $h$, which we presented in section \ref{sec:main3}.
We generated $500$ covariance matrices $\bC$ using the graphical model described in section \ref{sec:num_Laplacian}, with  $k=20$ and $p= 2,000$. Figure \ref{fig:accuracy6} displays $P[\log(h)]$ using our theoretical distribution $f_H(h)$ given by \eqref{eq:main_density} (dashed red curve) as a function of $\log(h)$. We also display as a solid black line the limiting scaling behavior,  $f_H(h)\varpropto h^{-2}$, given by \eqref{eq:density_largep}. Finally, we compare
these two probability density functions with the empirical distribution of $\log(h)$, shown as blue crosses.
We note for large $h$ that all distributions are parallel in this log-log plot, indicating that they have the same asymptotic power-law scaling.

%________________________________________________________________
\section{Discussion}\label{sec:discussion}
%________________________________________________________________
A central motivator for our research has been equation \eqref{eq:error}, which describes the limiting $n \rightarrow \infty$ expected sample error $ \|{\bu}_i-\bm{\tilde {u}}_i\|^2$ of a sample eigenvector $\tilde{\bu}_i$ for a covariance matrix drawn from a Wishart distribution. However, this equality only holds asymptotically. In this Discussion, we describe the conditions in which the approximation \eqref{eq:expectdelta} is expected to be accurate. That is, when is the sample size $n$ sufficiently large for given covariance matrix size $p$?

The standard approach to this problem usually involves a tail bound. Instead, we use here a simple argument that yields a lower bound that works very well in practice. Indeed, we provide a necessary (but not sufficient) lower bound on $n$ such that \eqref{eq:error} and \eqref{eq:main_approx} are valid. Since both ${\bu}_i$ and $\bm {\tilde {u}}_i$ are normalized and we assume ${\bu}_i \approx \bm{\tilde {u}}_i$, we have
\begin{equation}
  \|{\bu}_i-\bm{\tilde {u}}_i\|^2 = 2[1-  \langle{\bu}_i,\bm {\tilde u_i}\rangle] \le 2
  \label{eq:bound}. 
\end{equation}
Under the approximation $ \|{\bu}_i-\bm{\tilde {u}}_i\|^2 \approx h_i/n$ given by \eqref{eq:expectdelta}, it follows that
\begin{equation}
 n\ge   h_i/2. \label{eq:bound_heuristic}
\end{equation}

We now provide numerical support for this bound using the graphical model described in section \ref{sec:num_Laplacian}. We first generate a $k$-regular graph with $p$ vertices and compute the unnormalized Laplacian matrix, $\bC$, which we treat as a covariance matrix. Let $\{\bu_i\}_{i=1}^p$ be the eigenvectors of $\bC$. In order to study the convergence of the empirical eigenvectors, we generate 100 random matrices $\widetilde{\bC}$ from the Wishart distribution $W(\bC,n)$.  For each random realization $\widetilde{\bC}$, we calculate its eigenvectors $\{\tu_i\}_{i=1}^p$ and compute the residual error $ \bu_i - \tu_i$ between the sample eigenvectors and the population eigenvectors.

\begin{figure}[h] 
\centering
\includegraphics[width=0.45\textwidth]{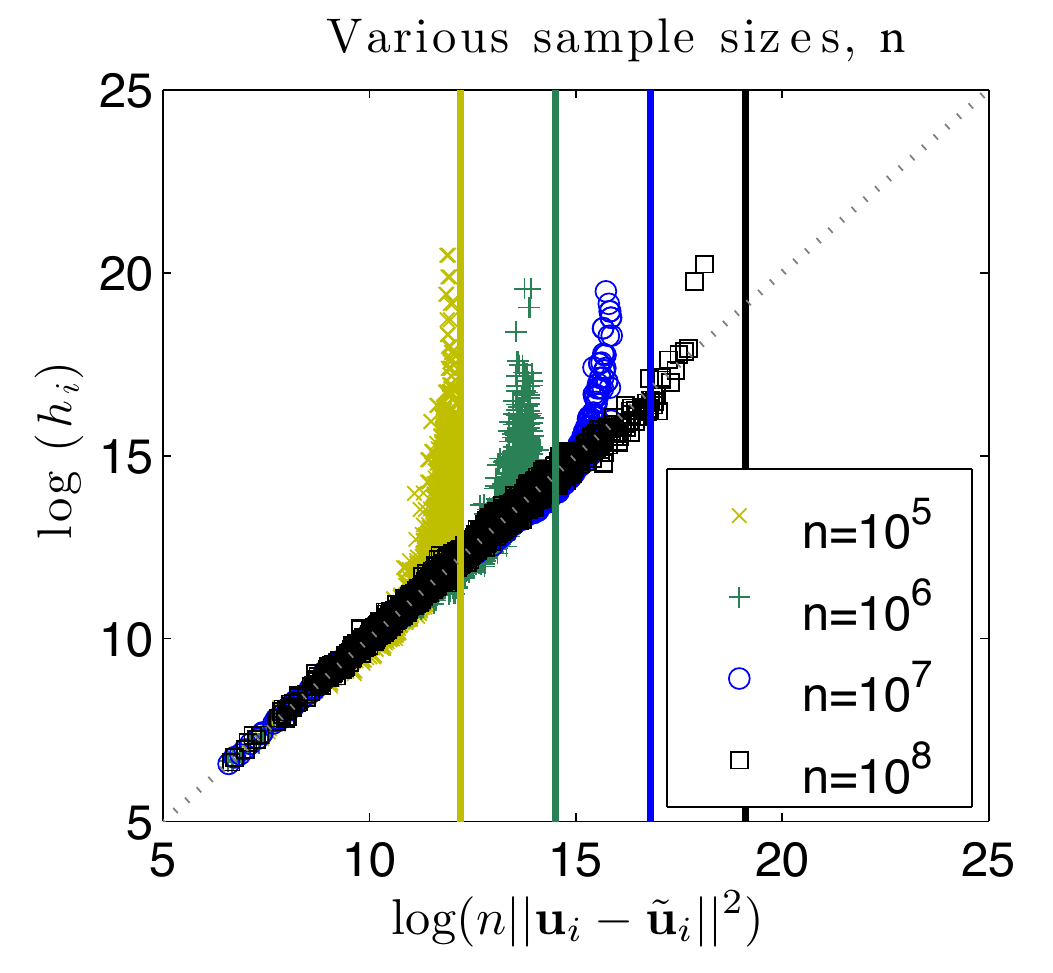}
\includegraphics[width=0.45\textwidth]{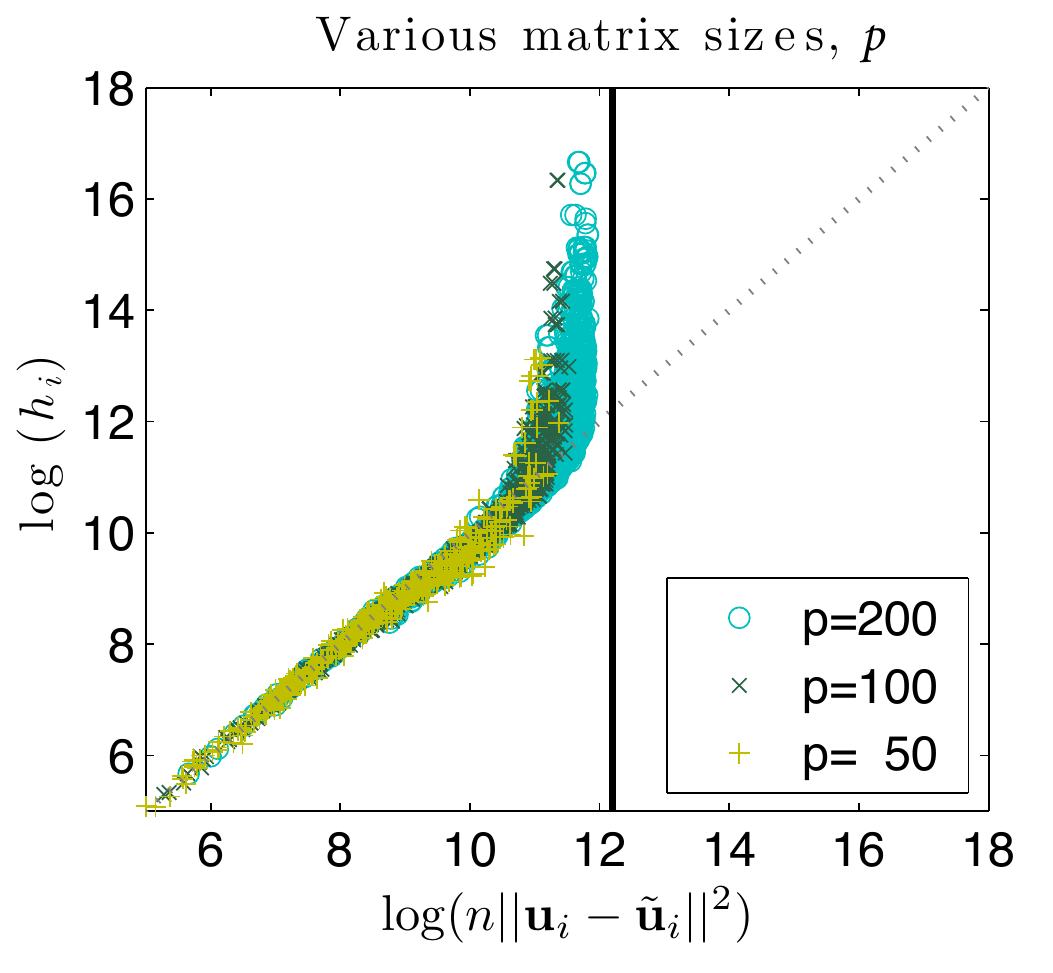}
\caption{
True {values} of $h_i$ given by \eqref{eq:hi} versus residual eigenvector error 
$n {\|\bu_i - \tu_i\|^2}$ across 100 Wishart-distributed sample covariance matrices $\widetilde{\bC}$ with expectation ${\bC}$ with (left) $p=200$ and various $n$,  and (right) $n=10^5$ and various $p$.
{The vertical lines indicate  $h_i\le 2n$ is a necessary  (but not sufficient) condition for accuracy of the approximation $\E{n\|\bu_i - \tu_i\|^2} \approx h_i/n$  given by \eqref{eq:expectdelta}.}
}
\label{fig:accuracy}
\end{figure}

Figure \ref{fig:accuracy}-left displays $\log(h_i)$ as a function of $\log(n { {\|\bu_i - \tu_i \|^2}})$ for each
random realization of a Wishart matrix $\widetilde{\bC}$ for $p=200$, $k=5$, and several choices of $n$. For each value of $n$, we plot the bound given by \eqref{eq:bound_heuristic}, $\log(2n)$, as a vertical solid line. 
Figure \ref{fig:accuracy}-right displays a scatterplot of $\log(h_i)$ as a function of $\log(n  {\|\bu_i - \tu_i \|^2} )$ for $k=5$, $n=10^5$, and several values of $p$. We also plot  $\log(2n)$ as a vertical solid line. 
Both panels  illustrate \eqref{eq:bound_heuristic} as a useful bound for considering when the approximation $h_i \approx   \E{n\|\bu_i - \tu_i \|^2}$ given by \eqref{eq:expectdelta} will be valid. Specifically, we require $ \E{\|\bu_i - \tu_i \|^2} < 2$ and observe considerable discrepancy as $\E{\|\bu_i - \tu_i \|^2} \to 2$.

\begin{figure}[t] 
\centering
\includegraphics[width=0.45\textwidth]{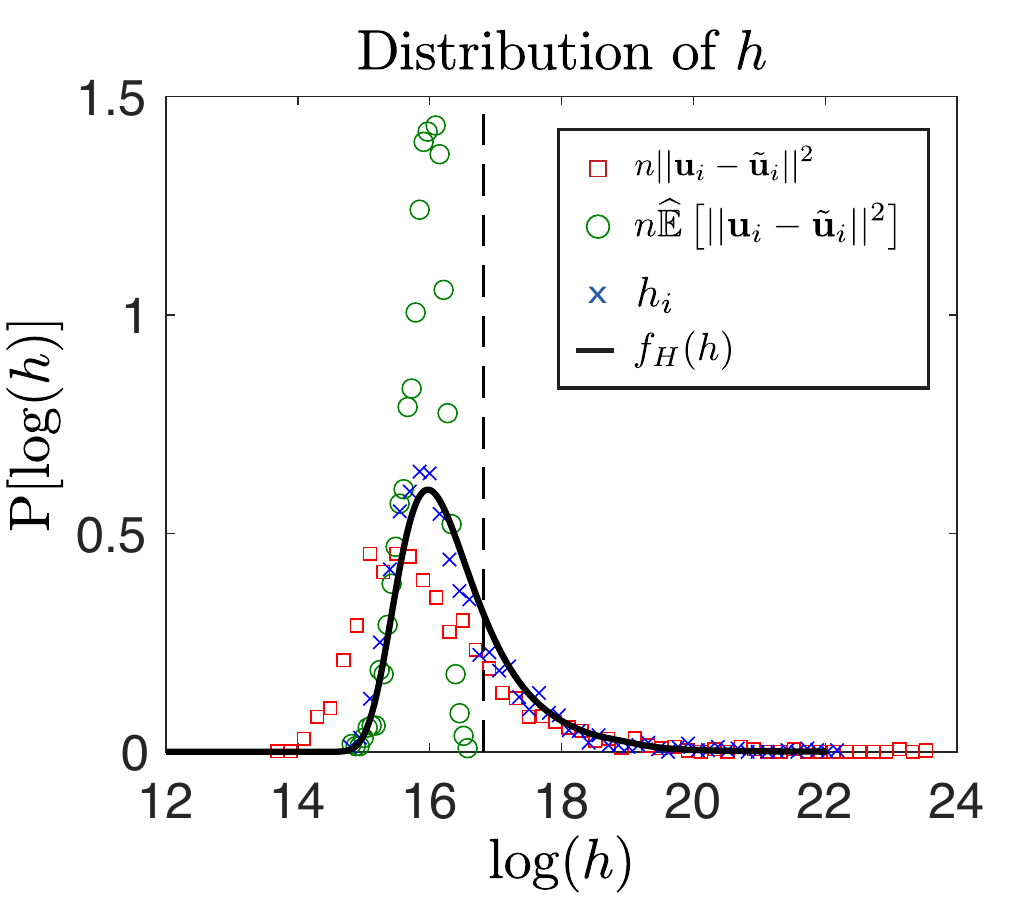}
\includegraphics[width=0.45\textwidth]{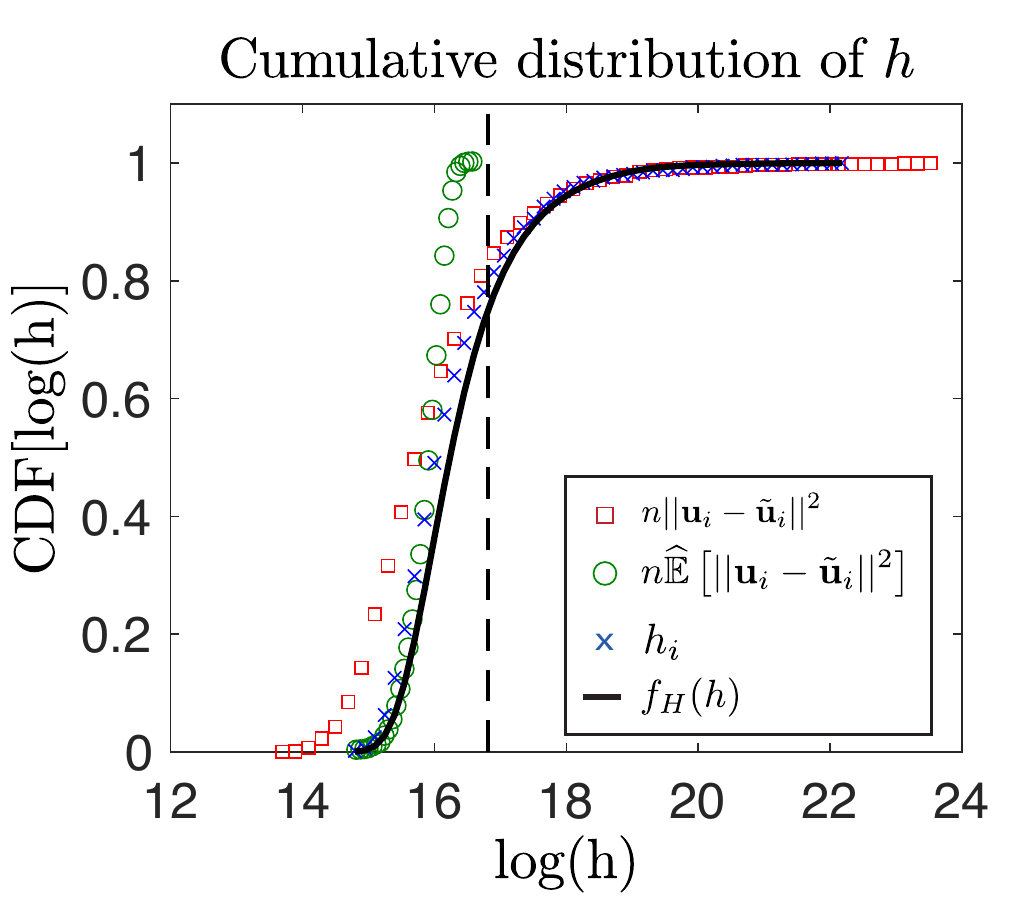}
\caption{
Discrepancy between theoretical distribution $f_H(h)$ of $h$ given by \eqref{eq:main_density}, 
%across a covariance matrix ensemble 
and a distribution of bootstrap estimates, $\wE{n\|\bu_i - \tu_i\|^2}$, due to error for  \eqref{eq:error} and  \eqref{eq:expectdelta}.
We plot results for two sources of error.
(green circles, ${\green \circ}$) Sample size $n$ is too small and does not satisfy bound \eqref{eq:bound_heuristic}, which is shown by the vertical dashed lines.
(red squares, $\red \Box$): Insufficiently many (particularly, $R=1$) samples are used to provide a reliable bootstrap estimate  $\wE{n\|\bu_i - \tu_i\|^2}$ to  $\E{n\|\bu_i - \tu_i\|^2}$. (See text for details.)
}
\label{fig:accuracy5}
\end{figure}

We conclude by exploring the effect of inaccuracy for \eqref{eq:expectdelta} on the distribution of the bootstrap estimates, $n\wE{\|\bu_i - \tu_i\|^2}$, of the mean residual error across the covariance matrix ensemble. We identify two sources of discrepancy: (i) choosing $n$ too small so that the bound \eqref{eq:bound_heuristic} is violated, and (ii) using insufficiently many samples from the Wishart distribution $W(\bC,n)$ to provide an accurate bootstrap estimate, $n\wE{\|\bu_i - \tu_i\|^2}$.
That is, the bootstrap estimate $\wE{\|\bu_i - \tu_i\|^2}$ is only a reliable estimate of $\E{\|\bu_i - \tu_i\|^2}$ if we generate enough random samples $\widetilde{\bC}$ from the Wishart distribution $W(\bC,n)$.

We highlight these two sources of discrepancy with a numerical experiment similar to the one described in section \ref{sec:num_4}. We generated 50 covariance matrices, $\bC$, with $k=20$ and $p=1,000$, as explained in section \ref{sec:num_Laplacian}. For each $\bC$,  we generated $R$ random realizations $\widetilde{\bC}$ from the Wishart distribution $W(\bC,n)$. For each $\widetilde{\bC}$, we calculated its eigenvectors $\{\tu_i\}_{i=1}^p$ and computed the residual error $\bu_i - \tu_i$ between the sample eigenvectors and the population eigenvectors. We then computed a bootstrap estimate, $\wE{\|\bu_i - \tu_i \|^2}$, of the mean sample error for each $\bC$ using the $R$ realizations of $\widetilde{\bC}$.

In Figure \ref{fig:accuracy5}-left and Figure \ref{fig:accuracy5}-right, we plot as solid black curves the probability distribution  $f_H(h)$ given by \eqref{eq:main_density}, and its corresponding cumulative distribution, respectively. We plot by blue crosses the empirical distribution of ${h}_i$, which we estimate using the 50 covariances $\bC$. Note that $f_H(h)$ accurately predicts the  observed distribution of $h_i$, since $p$ is sufficiently large. In addition, we plot the distribution of the bootstrap estimates $\wE{n\|\bu_i - \tu_i\|^2}$ of the mean sample error for  $R=10$ and $n=10^7$ (green circles) as well as $R=1$ and $n=10^{10}$ (red squares). Note that when $R=1$, the bootstrap estimate $\wE{\|\bu_i - \tu_i\|^2}$ is actually just the sample error ${\|\bu_i - \tu_i\|^2}$. 

Observe that both distributions disagree with $f_H(h)$ for different reasons: 
For $R=10$ and $n=10^7$, the distribution of $\wE{n\|\bu_i - \tu_i\|^2}$ is expected to differ because $n=10^7$ is too small and does not satisfy the bound given by \eqref{eq:bound_heuristic} (vertical dashed lines).  
On the other hand, for $R=1$ and $n=10^{10}$, the sample error does not provide a good bootstrap estimate for the mean sample error $\E{n\|\bu_i - \tu_i\|^2}$, which is the relevant quantity that is described in \eqref{eq:error} and \eqref{eq:expectdelta}. We observe in Figure \ref{fig:accuracy5} that using too few samples (i.e., small $R$) affects  the distribution of $\wE{n\|\bu_i - \tu_i\|^2}$ by shifting it toward small  values of $h$ (see red squares).

%________________________________________________________________
\section*{Acknowledgment}
The work of D. T. was partially supported by NSF Grant DMS-1127914 to the Statistical and Applied Mathematical Sciences Institute (SAMSI). Any opinions, findings, and conclusions or recommendations expressed in this material are those of the author(s) and do not necessarily reflect the views of the NSF.

%________________________________________________________________
\bibliographystyle{imamat}
\bibliography{my_bib}

\begin{thebibliography}{}

\bibitem[Abul-Magd \& Simbel, 1999]{abul1999wigner}
Abul-Magd, A. {\&} Simbel, M. (1999)  Wigner surmise for high-order level
  spacing distributions of chaotic systems. {\em Physical Review E},
  \textbf{60}(5), 5371.

\bibitem[Akemann et~al., 2010]{akemann2010universal}
Akemann, G., Fischmann, J. {\&} Vivo, P. (2010)  Universal correlations and
  power-law tails in financial covariance matrices. {\em Physica A: Statistical
  Mechanics and its Applications}, \textbf{389}(13), 2566--2579.

\bibitem[Anderson, 2003]{anderson03}
Anderson, T.~W. (2003) {\em An Introduction to Multivariate Statistical
  Analysis}.
John Wiley and Son, Inc., Hoboken, New Jersey, 3rd ed.

\bibitem[Bapat, 2010]{Bapat10}
Bapat, R.~B. (2010) {\em Graphs and Matrices}.
Springer.

\bibitem[Bassett et~al., 2011]{Bassett2011}
Bassett, D.~S., Wymbs, N.~F., Porter, M.~A., Mucha, P.~J., Carlson, J.~M. {\&}
  Grafton, S.~T. (2011)  Dynamic reconfiguration of human brain networks during
  learning. {\em Proceeding of the National Academy of Sciences},
  \textbf{108}(18), 7641--7646.

\bibitem[Benaych-Georges \& Nadakuditi, 2011]{benaych2011eigenvalues}
Benaych-Georges, F. {\&} Nadakuditi, R.~R. (2011)  The eigenvalues and
  eigenvectors of finite, low rank perturbations of large random matrices. {\em
  Advances in Mathematics}, \textbf{227}(1), 494--521.

\bibitem[Bickel \& Levina, 2008]{Bickel2008}
Bickel, P.~J. {\&} Levina, E. (2008)  Covariance regularization by
  thresholding. {\em The Annals of Statistics}, pages 2577--2604.

\bibitem[Bourgade et~al., 2014]{bourgade14}
Bourgade, P., Erd{\"o}s, L. {\&} Yau, H.-T. (2014)  Edge universality of beta
  ensembles. {\em Communications in Mathematical Physics}, \textbf{332}(1),
  261--353.

\bibitem[Brody, 1973]{brody1973statistical}
Brody, T. (1973)  A statistical measure for the repulsion of energy levels.
  {\em Lettere al Nuovo Cimento (1971-1985)}, \textbf{7}(12), 482--484.

\bibitem[Chung et~al., 2003]{Chung2003}
Chung, F., Lu, L. {\&} Vu, V. (2003)  Spectra of random graphs Spectra of
  random graphs with given expected degrees. {\em Proceeding of the National
  Academy of Sciences, USA}, \textbf{100}, 6313--6318.

\bibitem[Delvenne et~al., 2010]{Delvenne2010}
Delvenne, J.-C., Yaliraki, N., Sophia, N. {\&} Barahona, M. (2010)  Stability
  of graph communities across time scales. {\em Proceedings of the National
  Academy of Sciences}, \textbf{107}(29), 12755--12760.

\bibitem[Dumitriu \& Pal, 2012]{dumitriu2012sparse}
Dumitriu, I. {\&} Pal, S. (2012)  Sparse regular random graphs: {S}pectral
  density and eigenvectors. {\em The Annals of Probability}, \textbf{40}(5),
  2197--2235.

\bibitem[Ellegaard et~al., 1995]{ellegaard1995spectral}
Ellegaard, C., Guhr, T., Lindemann, K., Lorensen, H., Nyg{\aa}rd, J. {\&}
  Oxborrow, M. (1995)  Spectral statistics of acoustic resonances in aluminum
  blocks. {\em Physical review letters}, \textbf{75}(8), 1546.

\bibitem[Elton et~al., 2009]{Elton2009}
Elton, E.~J., Gruber, M.~J., Brown, S.~J. {\&} Goetzmann, W.~N. (2009) {\em
  Modern portfolio theory and investment analysis}.
John Wiley and Sons.

\bibitem[Gatti et~al., 2010]{Gatti2010}
Gatti, D.~M., Barry, W.~T., Nobel, A.~B., Rusyn, I. {\&} Wright, F.~A. (2010)
  Heading down the wrong pathway: on the influence of correlation within gene
  sets. {\em BMC Genomics}, \textbf{11}(1), 574.

\bibitem[Golub \& Loan, 2012]{Golub2012}
Golub, G.~H. {\&} Loan, C. F.~V. (2012) {\em Matrix Computations}, volume~3.
JHU Press.

\bibitem[Guhr et~al., 1998]{guhr1998random}
Guhr, T., M{\"u}ller-Groeling, A. {\&} Weidenm{\"u}ller, H.~A. (1998)
  Random-matrix theories in quantum physics: common concepts. {\em Physics
  Reports}, \textbf{299}(4-6), 189--425.

\bibitem[Hastie et~al., 2009]{Hastie2009}
Hastie, T., Tibshirani, R. {\&} Friedman, J. (2009) {\em The Elements of
  Statistical Learning}, volume~2.
New York: Springer.

\bibitem[Herman et~al., 2007]{Herman2007}
Herman, D., Ong, T.~T., Usaj, G., Mathur, H. {\&} Baranger, H.~U. (2007)  Level
  Spacings in Random Matrix Theory and Coulomb Blockade Peaks in Quantum Dots.
  {\em Physical Review B}, \textbf{76}, 195448.

\bibitem[I.~J.~Farkas \& Vicsek, 2001]{Farkas2001}
I.~J.~Farkas, I.~Der\'enyi, A.-L.~B. {\&} Vicsek, T. (2001)  Beyond the
  semicircle Beyond the semicircle law. {\em Physical Review E}, \textbf{64},
  026704.

\bibitem[K.-I.~Goh \& Kim, 2001]{Goh2001}
K.-I.~Goh, B.~K. {\&} Kim, D. (2001)  Spectra and eigenvectors of scale-free
  networks. {\em Physical Review E}, \textbf{64}, 051903.

\bibitem[Kuhn, 2008]{Kuhn2008}
Kuhn, R. (2008)  Spectra of sparse random matrices. {\em Journal of Physics A},
  \textbf{41}, 295002.

\bibitem[Lam \& Fan, 2009]{Lam2009}
Lam, C. {\&} Fan, J. (2009)  Sparsistency and rates of convergence in large
  covariance matrix estimation. {\em Annals of Statistics}, \textbf{37}(6B),
  4254.

\bibitem[Ledoit \& P{\'e}ch{\'e}, 2011]{ledoit2011eigenvectors}
Ledoit, O. {\&} P{\'e}ch{\'e}, S. (2011)  Eigenvectors of some large sample
  covariance matrix ensembles. {\em Probability Theory and Related Fields},
  \textbf{151}(1-2), 233--264.

\bibitem[Mantegna \& Stanley, 2000]{Mantegna2000}
Mantegna, R. {\&} Stanley, H. (2000) {\em An Introduction to Econophysics}.
Cambridge University Press: Cambridge, UK.

\bibitem[McKay, 1981]{McKay1981}
McKay, B. (1981)  The expected eigenvalue distribution of a large regular
  graph. {\em Journal of Linear Algebra and its Applications}, \textbf{40},
  203--216.

\bibitem[Mehta, 1991]{Mehta1991}
Mehta, M.~L. (1991) {\em Random Matrices, 3rd ed.}
Academic Press, New York.

\bibitem[Murphy et~al., 2017]{murphy2017anderson}
Murphy, N.~B., Cherkaev, E. {\&} Golden, K.~M. (2017)  Anderson transition for
  classical transport in composite materials. {\em Physical review letters},
  \textbf{118}(3), 036401.

\bibitem[Newman, 2003]{Newman2003}
Newman, M. E.~J. (2003)  The structure and function of complex networks. {\em
  SIAM Review}, \textbf{45}(2), 167--256.

\bibitem[Pastur et~al., 2011]{pastur11}
Pastur, L.~A., Shcherbina, M. {\&} Shcherbina, M. (2011) {\em Eigenvalue
  distribution of large random matrices}, volume 171.
American Mathematical Society Providence, RI.

\bibitem[Peixoto, 2013]{peixoto2013eigenvalue}
Peixoto, T.~P. (2013)  Eigenvalue spectra of modular networks. {\em Physical
  Review Letters}, \textbf{111}(9), 098701.

\bibitem[Pimpinelli et~al., 2005]{pimpinelli2005evolution}
Pimpinelli, A., Gebremariam, H. {\&} Einstein, T. (2005)  Evolution of
  terrace-width distributions on vicinal surfaces: Fokker-Planck derivation of
  the generalized Wigner surmise. {\em Physical review letters},
  \textbf{95}(24), 246101.

\bibitem[Plerou et~al., 2002]{plerou2002random}
Plerou, V., Gopikrishnan, P., Rosenow, B., Amaral, L. A.~N., Guhr, T. {\&}
  Stanley, H.~E. (2002)  Random matrix approach to cross correlations in
  financial data. {\em Physical Review E}, \textbf{65}(6), 066126.

\bibitem[S.~N.~Dorogovtsev \& Samukhin, 2003]{Dorogovtsev2003}
S.~N.~Dorogovtsev, A. V.~Goltsev, J. F. F.~M. {\&} Samukhin, A.~N. (2003)
  Spectra of complex networks. {\em Physical Review E}, \textbf{68}, 046109.

\bibitem[Schierenberg et~al., 2012]{schierenberg2012wigner}
Schierenberg, S., Bruckmann, F. {\&} Wettig, T. (2012)  Wigner surmise for
  mixed symmetry classes in random matrix theory. {\em Physical Review E},
  \textbf{85}(6), 061130.

\bibitem[Shklovskii et~al., 1993]{shklovskii1993statistics}
Shklovskii, B., Shapiro, B., Sears, B., Lambrianides, P. {\&} Shore, H. (1993)
  Statistics of spectra of disordered systems near the metal-insulator
  transition. {\em Physical Review B}, \textbf{47}(17), 11487.

\bibitem[Taylor et~al., 2017]{taylor2017super}
Taylor, D., Caceres, R.~S. {\&} Mucha, P.~J. (2017)  Super-resolution community
  detection for layer-aggregated multilayer networks. {\em Physical Review X},
  \textbf{7}(3), 031056.

\bibitem[Taylor et~al., 2016]{taylor2016enhanced}
Taylor, D., Shai, S., Stanley, N. {\&} Mucha, P.~J. (2016)  Enhanced
  detectability of community structure in multilayer networks through layer
  aggregation. {\em Physical Review Letters}, \textbf{116}(22), 228301.

\bibitem[Volkov et~al., 2009]{Volkov2009}
Volkov, I., Banavar, J.~R., Hubbell, S.~P. {\&} Maritan, A. (2009)  Inferring
  species interactions in tropical forests. {\em Proceedings of the National
  Academy of Sciences}, \textbf{106}(33), 13854--13859.

\bibitem[Weigt et~al., 2009]{Weigt2009}
Weigt, M., White, R.~A., Szurmant, H., Hoch, J.~A. {\&} Hwa, T. (2009)
  Identification of direct residue contacts in protein--protein interaction by
  message passing. {\em Proceedings of the National Academy of Sciences},
  \textbf{106}(1), 67--72.

\bibitem[Wigner, 1958]{wigner1958distribution}
Wigner, E.~P. (1958)  On the distribution of the roots of certain symmetric
  matrices. {\em Annals of Mathematics}, pages 325--327.

\bibitem[Wigner, 1993]{wigner1993class}
Wigner, E.~P. (1993)  On a class of analytic functions from the quantum theory
  of collisions. In {\em The Collected Works of Eugene Paul Wigner}, pages
  409--440. Springer.

\bibitem[Zhang et~al., 2014]{Zhang2004}
Zhang, X., Nadakuditi, R.~R. {\&} Newman, M. E.~J. (2014)  Spectra of random
  graphs with community structure and arbitrary degrees. {\em Physical Review
  E}, \textbf{89}, 042816.

\end{thebibliography}
%________________________________________________________________

\section*{Appendices}

%________________________________________________________________
\begin{appendices}
%________________________________________________________________

%________________________________________________________________
\section{Derivation of main result 1\label{sec:proof_main1}}
%________________________________________________________________

In this Appendix, we approximate $h_i$ in \eqref{eq:error} in terms of nearest neighbor eigenvalue
gaps. By doing so, we will be able to exploit the knowledge of the $p\to\infty$ limiting distribution of the eigenvalues.
We begin by dividing the summation into two parts so that
\begin{equation}
   h_i =  h_i^-+ h_i^+,
\end{equation}
with
\begin{eqnarray}
   h_i^-&=& \sum_{j=1}^{i-1} \frac{\lambda_i\lambda_j}{(\lambda_i-\lambda_j)^2},\label{eq:h1}\\
   h_i^+ &=& \sum_{j=i+1}^{p} \frac{\lambda_i\lambda_j}{(\lambda_i-\lambda_j)^2}  .\label{eq:h2}
\end{eqnarray}
Our numerical experiments show that typically the nearest neighbor terms dominate the others. Taking this into account, we isolate the first spacing, and rewrite $h_i^{\pm}$ as
\begin{eqnarray}
   h_i^{-}&=&  \frac{\lambda_i\lambda_{i-1}}{(\lambda_i-\lambda_{i-1})^2} + 
   \sum_{j=1}^{i-2} \frac{\lambda_i\lambda_ j}{(\lambda_i-\lambda_ j)^2} \label{eq:removed1},\\
   h_i^{+}&=&  \frac{\lambda_i\lambda_{i+1}}{(\lambda_i-\lambda_{i+1})^2} + 
   \sum_{j=i+2}^p \frac{\lambda_i\lambda_j}{(\lambda_i-\lambda_j)^2}.\label{eq:removed2}
\end{eqnarray}

We study the large $p$ behavior of \eqref{eq:removed1} and \eqref{eq:removed2} by separately considering the nearest-neighbor terms and the summations. In particular, we will obtain approximations that rely only on the right and left nearest-neighbor eigenvalue gaps,
\begin{eqnarray}
   s_i^{\pm}=|\lambda_i- \lambda_{i\pm 1}|.
\end{eqnarray}
We first consider the isolated terms, 
\begin{eqnarray}
  \frac{\lambda_i\lambda_{i\pm1}}{(\lambda_i-\lambda_{i\pm1})^2} 
  &=& \frac{\lambda_i(\lambda_i \pm s_i^{\pm} )}{(s_i^{\pm})^2} \\
  &=&\frac{\lambda_i^2}{(s_i^{\pm})^2} \left[ 1+\mathcal{O}(s_i^{\pm})\right] . 
\end{eqnarray}
Using that $s_i^{\pm}\to 0 $ as $p\to\infty$ (which is established by assumption \ref{asssume3} and  convergences, in expectation, with rate $s_i^{\pm}=\mathcal{O}(1/p)$), we find the asymptotic estimate
\begin{eqnarray}
   \frac{\lambda_i\lambda_{i\pm1}}{(\lambda_i-\lambda_{i\pm1})^2}  \to \frac{\lambda_i ^2}{(s_i^{\pm})^2}  .
\end{eqnarray}

We now turn our attention to the summations, which we will approximate using the limiting $p\rightarrow \infty$ spectral density $\rho(\lambda)$ of the normalized empirical counting measure of the eigenvalues. 
More precisely, consider a sequence of size-$p$ symmetric covariance matrices, each having eigenvalues $\{\lambda_j\}$ for $j\in\{1,\dots,p\}$. We define for each matrix the empirical spectral density 
\begin{equation}
  \rho_p(\lambda) = p^{-1}\sum_j \delta({\lambda_j }),\label{eq:density}
\end{equation}
where $\delta(\lambda)$ is the Dirac delta function and $\lambda\in\mathbb{R}$. We assume the covariance matrices are drawn from an ensemble such that the sequence $\{\rho_p(\lambda)\}$ weakly converges, implying that
\begin{equation}
  \int_{-\infty}^\infty \rho_p(\lambda)f(\lambda) d\lambda \to \int_{-\infty}^\infty \rho(\lambda)f(\lambda) d\lambda  
  \label{eq:weak_converg}
\end{equation}
as $p\to\infty$ for any continuous and bounded function $f(\lambda)$. We assume $\rho(\lambda)$ is continuous, bounded, has compact support (denoted $\text{supp}(\rho)$), and is differentiable on $\text{supp}(\rho)$. For notational convenience, we assume $\text{supp}(\rho)=(\alpha,\beta)$ for some $\alpha,\beta\in\mathbb{R}$, allowing us to replace the limits of integration in \eqref{eq:weak_converg} by $(\alpha,\beta)$. However, our analysis can be easily extended to unions of such intervals.

We begin be rewriting the summations in \eqref{eq:removed1} and \eqref{eq:removed2} as the integration of function 
\begin{eqnarray}
  f_{\lambda_i}(\lambda)&=\frac{\lambda_i\lambda}{(\lambda_i-\lambda)^2} \label{eq:def_fi}.
\end{eqnarray}
with probability measure $\rho_p(\lambda)$ given by \eqref{eq:density},
\begin{eqnarray}
  \frac{1}{p}\sum_{j=1}^{i-2} \frac{\lambda_i\lambda_ j}{(\lambda_i-\lambda_ j)^2} &=& 
  \int_{\alpha}^{\lambda_{i-1}} \rho_p(\lambda)f_{\lambda_i}(\lambda) d\lambda, \label{eq:int1}   \\
  \frac{1}{p}\sum_{j=i+2}^p \frac{\lambda_i\lambda_j}{(\lambda_i-\lambda_j)^2}&=& 
  \int_{\lambda_{i+1}}^\beta \rho_p(\lambda)f_{\lambda_i}(\lambda) d\lambda .\label{eq:int2}
\end{eqnarray}
Because $f_{\lambda_i}(\lambda)$ is unbounded at the singularity $\lambda=\lambda_i$, \eqref{eq:weak_converg} does not describe the behavior of integral $ \int_{\alpha}^\beta \rho_p(\lambda)f_{\lambda_i}(\lambda)d\lambda$, which we find to diverge with $p$ for any  $\lambda_i\in \text{supp}(\rho)$. Fortunately, \eqref{eq:int1} and \eqref{eq:int2} do not require integration across the singularity at $\lambda=\lambda_i$; however, the limits of integration, i.e., $\lambda_{i-1}$ in \eqref{eq:int1} and $\lambda_{i+1}$ in \eqref{eq:int2}, depend on $p$ (and converge to the singularity at $\lambda_i$). Thus, \eqref{eq:weak_converg} is also not directly applicable to \eqref{eq:int1} and \eqref{eq:int2}.

To proceed, we restrict our attention to \eqref{eq:int1} since analogous results can be obtained for \eqref{eq:int2}. We consider, for the moment, \eqref{eq:int1} with fixed upper limit $\lambda_i-\epsilon$ for $\epsilon>0$ and $\epsilon\approx0$. Equation \eqref{eq:weak_converg} implies the $p\to\infty$ limit 
\begin{equation}
\int_{\alpha }^{\lambda_i -\epsilon} f_{\lambda_i}(\lambda)\rho_p(\lambda)d\lambda \to \int_{\alpha }^{\lambda_i -\epsilon} f_{\lambda_i}(\lambda)\rho(\lambda)d\lambda  .\label{eq:lim_lower}
\end{equation}
We now study how the right-hand side of \eqref{eq:lim_lower} scales with $\epsilon$. 
Using that both $\rho(\lambda)$ and $f_{\lambda_i}(\lambda)$ are differentiable for $\lambda\in\text{supp}(\rho)\setminus \{\lambda_i\}$,  we implement integration by parts, treating the numerator and denominator separately, to obtain
\begin{eqnarray}
   \int_{\alpha }^{\lambda_i -\epsilon} f_{\lambda_i}(\lambda)\rho(\lambda)d\lambda 
   &=& \lambda_i\frac{(\lambda_i-\epsilon)\rho(\lambda_i-\epsilon)}{\epsilon}  - 
\lambda_i\int_{\alpha }^{\lambda_i -\epsilon}  \frac{\rho(\lambda) + \lambda\rho'(\lambda)}{\lambda_i-\lambda} d\lambda .
   \label{eq:int_byparts}
\end{eqnarray}
The first term in the right hand side of \eqref{eq:int_byparts} has the $\epsilon\to0$ asymptotic estimate
\begin{eqnarray}
\lambda_i\frac{(\lambda_i-\epsilon)\rho(\lambda_i-\epsilon)}{\epsilon}  \to \frac{\lambda_i^2\rho(\lambda_i)}{\epsilon} .
\label{eq:first_term}
\end{eqnarray}
The second term on the right-hand side of \eqref{eq:int_byparts} is bounded as
\begin{alignat}{2}
   \left| \lambda_i\int_{\alpha}^{\lambda_{i}-\epsilon} \frac{\left[\rho(\lambda) +\lambda\rho'(\lambda)\right]}{\lambda_i-\lambda} d   
 \lambda \right|
  & \leq \lambda_i\left[\sup_{\lambda \in (\alpha, \lambda_{i}-\epsilon]} |\rho(\lambda) +\lambda\rho'(\lambda)|\right]
\int_{\alpha}^{\lambda_{i}-\epsilon} \frac{1}{|\lambda_i-\lambda|} d\lambda \nonumber\\
   & = \lambda_i\left[\sup_{\lambda \in (\alpha, \lambda_{i}-\epsilon]} |\rho(\lambda) +\lambda\rho'(\lambda)\right]
\ln \left(\frac{\lambda_i-\alpha}{\epsilon}\right)
   \label{eq:sup_bound}
\end{alignat}
It follows that the second term in the right-hand side of \eqref{eq:int_byparts} has scaling $\mathcal{O}(\text{ln}(1/\epsilon))$ and is dominated in the limit $\epsilon\to0$ by the first term, which is $\mathcal{O}(1/\epsilon)$. We combine \eqref{eq:first_term} and \eqref{eq:sup_bound}  to obtain the $\epsilon\to0$ asymptotic estimate
\begin{equation}
   \int_{\alpha }^{\lambda_i -\epsilon} f_{\lambda_i}(\lambda)\rho(\lambda)d\lambda \to 
   \frac{\lambda_i^2\rho(\lambda_i)}{\epsilon}  .
   \label{eq:eps_limit}
\end{equation}

We finally note that in the case where $\rho'(\lambda)$ is unbounded, it is straightforward to separate the integral on the left-hand side of \eqref{eq:sup_bound} into two domains, one containing all values $\lambda$ where $\rho'(\lambda)$ is unbounded and the second domain having upper limit $\lambda_i-\epsilon$. The first integral will converge to zero due to \eqref{eq:weak_converg}, whereas the second satisfies the bound given by \eqref{eq:sup_bound}, implying that the integral term in \eqref{eq:int_byparts} is $\mathcal{O}(\ln(1/\epsilon))$ provided that $\rho(\lambda)$ is differentiable in a small neighborhood containing $\lambda_i$.

We study the $p\to\infty$ limiting behavior for the right-hand side of \eqref{eq:int1} by considering the following identity,
\begin{align}
 \int_{\alpha}^{\lambda_{i-1}} f_{\lambda_i}(\lambda) \rho_p(\lambda) d\lambda 
 &=  \int_{\alpha}^{\lambda_{i}-s_i^-} f_{\lambda_i}(\lambda) \rho(\lambda) d\lambda +
 \int_{\alpha}^{\lambda_{i}-s_i^-}f_{\lambda_i}(\lambda) \left[\rho_p(\lambda)-\rho(\lambda) \right]  d\lambda .
 \label{eq:assumption_weak}
 \end{align}
The first term on the right-hand side grows linearly with $p$, which is straightforward to show by setting  $\epsilon=s_i^-$ in \eqref{eq:eps_limit} and using that $s_i^-=\mathcal{O}(1/p)$. Turning our attention to the second term on the right-hand side of \eqref{eq:assumption_weak}, recall that it would converge to zero if the upper limit of integration was fixed. However, $\lambda_{i}-s_i^-$ limits to $\lambda_i$ and the $p\to\infty$ behavior of the second term therefore depends on the rate of weak convergence for $\rho_p(\lambda) \to \rho(\lambda)$.  We assume that this term scales sublinearly with $p$ and is therefore dominated by the first term on the right-hand side of \eqref{eq:assumption_weak}. Under this assumption (and by conducting a similar analysis for \eqref{eq:removed2}), we obtain the asymptotic estimates
\begin{align}
p \int_{\alpha}^{\lambda_{i-1}} f_{\lambda_i}(\lambda) p \rho_p(\lambda) d\lambda 
& \to \frac{\lambda_i^2p\rho(\lambda_i)}{s_i^-}  , \label{eq:sym_two} \\
p \int_{\lambda_{i+1}}^\beta f_{\lambda_i}(\lambda) \rho_p(\lambda) d\lambda 
& \to \frac{\lambda_i^2p\rho(\lambda_i)}{s_i^+}  . \label{eq:sym_two2}
 \end{align}

%\todo{
%Note that this assumption is equivalent to
%\begin{align}\label{eq:ad}
%	 p^{-1}\int_\alpha^{\lambda_i-s^-} f_{\lambda_i}(\lambda) \left[\rho_p(\lambda)-\rho(\lambda) \right]  d\lambda 
%	 &\to  0
%	 %f_{\lambda_i}(\lambda) \left[ \int \rho_p(\lambda)-\rho(\lambda) d\lambda \right]  -  \int df_{\lambda_i}(\lambda) \int \left[\rho_p(\lambda)-\rho(\lambda) \right]  d\lambda \nonumber\\
%	 %\int f_{\lambda_i}(\lambda) \rho_p(\lambda)   d\lambda 
%	% - \int f_{\lambda_i}(\lambda)  \rho(\lambda)   d\lambda  \nonumber\\
%	%&= \int f_{\lambda_i}(\lambda) \rho_p(\lambda)   d\lambda 
%	% - \int f_{\lambda_i}(\lambda)  \rho(\lambda)   d\lambda  \nonumber\\
% \end{align}
% as $p\to \infty$ and $s^-\to 0 $ with $\mathbb{E}[s^-] = \mathcal{O}(p^{-1})$.
%}

In summary, we combine \eqref{eq:sym_two}, \eqref{eq:sym_two2} and \eqref{eq:first_term} to obtain the asymptotic  large $p$ approximation,
\begin{equation}
h_i^{\pm} \approx  \frac{\lambda^2_i}{(s_i^\pm)^2}  + \frac{p \rho(\lambda_i) \lambda^2_i}{s_i^\pm},
\label{eq:h_approx} 
\end{equation}
which gives the approximation \eqref{eq:main_approx}. We stress that this approximation assumes a sufficiently high rate of weak convergence for the spectral density so that the second term on the right-hand side of \eqref{eq:assumption_weak} is sublinear.

%________________________________________________________________
\section{Derivation of main result  2} \label{sec:estimate2}
%________________________________________________________________

In this section, we take a different perspective, and consider $h_i$, defined by \eqref{eq:hi}, to be
the realization of a random variable that is a function of the corresponding family of random
covariance matrices. Using the approximation $\hat{h}$ of $h$ provided by \eqref{eq:main_approx}, we
derive an estimate for the probability distribution, $P(h)$ of $h$. Let us denote by
$H$  the random variable for which $h_i$ is a realization.

Our goal is to remove the dependency on the random variables $s^+$ and $s^-$ in
\eqref{eq:main_approx}, so that $\hat{h}$ becomes a function of only $\lambda$, which is distributed
according the density $\rho(\lambda)$. The only missing ingredients are the probability
distributions of $s^+$ and $s^-$. We note that these two random variables are correlated, and thus
our line of attack involves using an approximation to the joint probability for the eigenvalue gaps,
$J(s^-,s^+)$, and derive an expression for the limiting probability density of the approximation
$\hat{h}$. In this section, we keep the discussion general, and derive an expression that is valid
for all $\rho(\lambda)$.

We assume that the joint probability distribution $J(s^-,s^+)$ of the left and right gaps around
each eigenvalue $\lambda$ can be approximated by \eqref{eq:joint}, which is reproduced below for ease
of presentation,
\begin{equation*}
J(s^- ,s^+) \approx 
 \frac{3^7\left[p \rho(\lambda)\right]^5}{32 \pi^3}
 \left[s^+s^- (s^+   +   s^- )\right] 
\exp\left(-\frac{\left[ 3p\rho(\lambda) \right]^2 }{4\pi}
 \left[(s^+ )^2 +(s^-)^2+s^+  s^- \right] \right). 
\end{equation*}
The expression \eqref{eq:joint} was derived in \cite{Herman2007} using $3 \times 3$ matrices from the
Gaussian Orthogonal Ensembles (GOE). As suggested by our numerical simulations (see Figure~\ref{joint_fig}), \eqref{eq:joint}
provides a good approximation for the covariance matrices that we study.

\begin{figure}[t] 
\centering 
\includegraphics[width=0.45\textwidth]{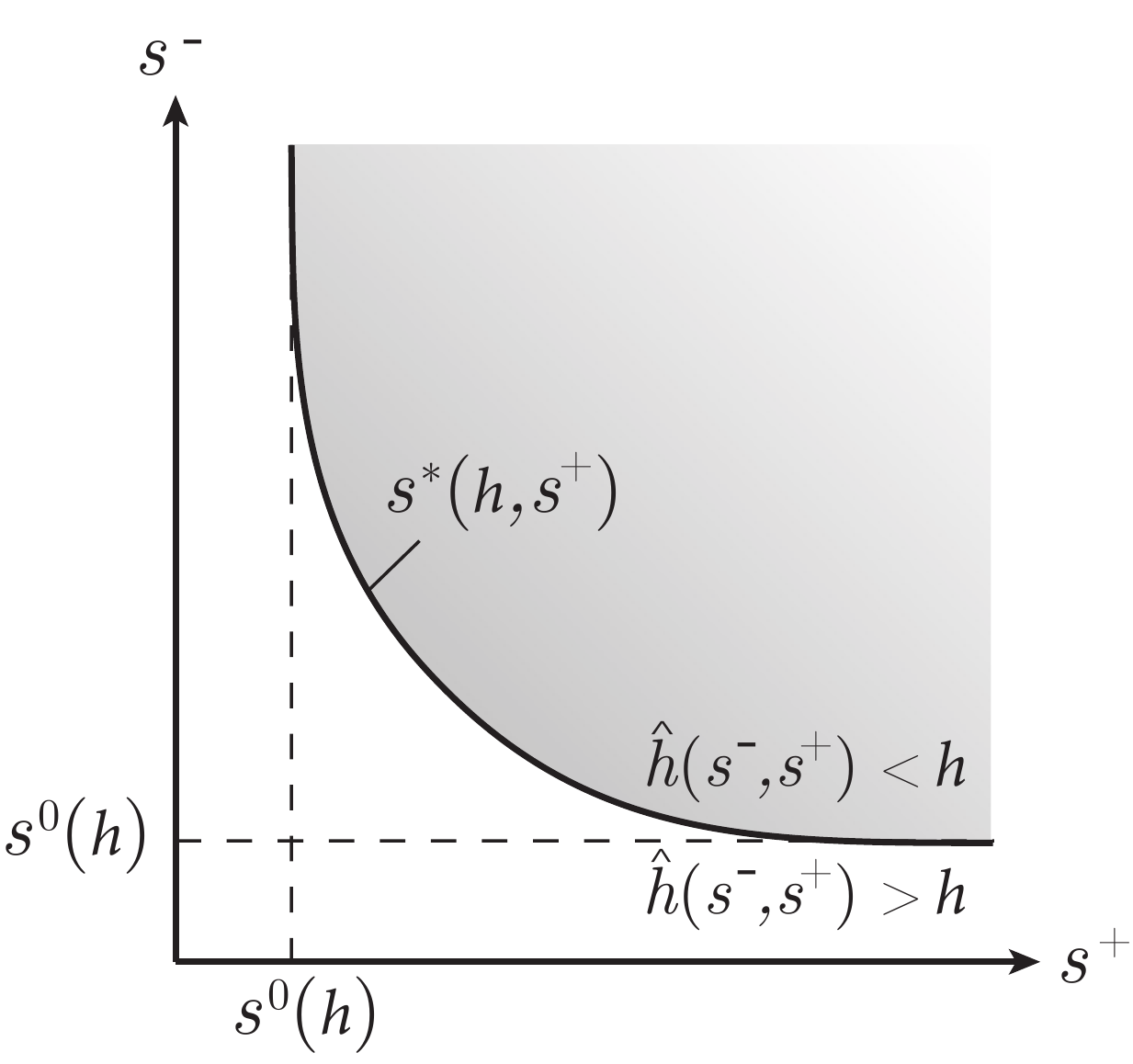}
\caption{
The cumulative distribution $F(h)$ for the random variable $H$ is shown as the integral
of $J(s^-,s^+)$ over region $\mathcal{S}$ given by \eqref{eq:S} (shaded
region). This region corresponds to $s^+\in(s^0(h),\infty)$ and $s^-\in\left(s^*(h,s^+),\infty\right)$, where $s^0(h)$ is found so that $h^+(s^+)<h$ for $s^+>s^0(h)$ and $s^*(h,s^+)$ is found so that $\hat{h}(s^-,s^+)<h$ for $s^->s^*(h,s^+)$. 
}
\label{fig:int_reg}
\end{figure}

To derive the distribution $P(h)$ of $H$, we first consider the cumulative distribution
\begin{equation}
   F(h) \stackrel{\Delta}{=} P(H < h).
\end{equation} 
Given an eigenvalue $\lambda$, we can find all the pairs of gaps $s^-$ and $s^+$,
such that $\hat{h}$ in \eqref{eq:main_approx} is less than $h$. Let
\begin{equation}
\mathcal{S}  \stackrel{\Delta}{=} \{(s^-,s^+):~\hat{h}(s^-,s^+) < h\} \label{eq:S}
\end{equation}
 be this set. We then proceed to compute the measure of
$\mathcal{S} $ using the joint probability density function defined above, 
\begin{equation}
F(h)  = \int_{\mathcal{S} }
J(s^-,s^+)\;ds^- ds^+.
\end{equation}
It turns out that we can describe analytically the set $\mathcal{S} $ (see
Figure~\ref{fig:int_reg}). For a given value of $h$, $\hat{h}(s^-,s^+)<h$ implies that both $h^+<h$
and $h^-<h$, where $h^\pm$ is given by \eqref{eq:h_approx} (with the subscript omitted), and
therefore the region of integration has the lower bounds $s^->s^0(h)$ and $s^+>s^0(h)$, where
$s^0(h)$ is given by
\begin{align} 
s^0(h) &= \frac{\lambda^2 p \rho(\lambda)}{2h} +
  \sqrt{\frac{\lambda^2}{h} +\left(\frac{\lambda^2 p \rho(\lambda)}{2h}\right)^2 }\nonumber\\
&= \frac{\lambda^2 p \rho(\lambda)}{2h}
\left(
1 + \sqrt{1 + \frac{4h}{\left[\lambda p \rho(\lambda)\right]^2}}
~~\right),
\end{align}
which follows directly from solving \eqref{eq:h_approx} for $s^\pm$ with $h^{\pm} = h$. We
therefore integrate $s^+$ over the range $(s^0(h),\infty)$.  For given values $h$ and $s^+$,
requiring that $\hat{h}(s^-,s^+)>h$ implies that $s^->s^*(s^+,h)$, where $s^*(h,s^+)$ is found by
substituting $h\mapsto \hat{h}(s^-,s^+)$ in \eqref{eq:main_approx} and solving for the positive root of $s^-$,
\begin{equation}
  s^*(h,s^+) = \lambda^2 p \rho(\lambda)  
  \frac{
    1 +  \sqrt{ 1 + \frac{\displaystyle 4}{\left[\displaystyle \lambda p \rho(\lambda)\right]^2}
       \left(
         h - \frac{\displaystyle \lambda^2}{\displaystyle (s^+)^2} - \frac{\displaystyle\lambda p \rho(\lambda) }{\displaystyle s^+} 
       \right)
    }
  }
  {
    2\left(h - \frac{\displaystyle \lambda^2}{\displaystyle (s^+)^2} -
      \frac{\displaystyle \lambda p \rho(\lambda) }{\displaystyle s^+}
    \right)
  }.
\end{equation}
We therefore integrate $s^-$ over the range $(s^*(h,s^+),\infty)$,
\begin{equation}
   F(h) = \int_{s^0(h)}^{\infty} \int_{{s^*}(h,s^+)}^{\infty} J(s^-,s^+)ds^- ds^+.\label{eq:cumulative}
\end{equation}

To obtain an estimate for the distribution of $h$, $f_H(h)$, we differentiate \eqref{eq:cumulative} with respect to $h$ to obtain
\begin{align}
f_H(h) & =\frac{\partial}{\partial h}  \int_{s^0(h)}^{\infty} \int_{ s^*(h,s^+)}^{\infty}
        J(s^-,s^+) ds^- ds^+ \label{eq:p_h2}\\ 
&=  - \frac{\partial s^0}{\partial h}(h) \int_{ s^*(h,s^0(h))}^{\infty}  J(s^-,s^0(h)) ds^- 
+ \int_{s^0(h)}^{\infty} \frac{\partial}{\partial h} \left[  \int_{ s^*(h,s^+)}^{\infty} 
  \mspace{-2mu} J(s^-,s^+)  ds^-  \right] ds^+\nonumber\\ 
&=  - \int_{s^0(h)}^{\infty}   J\left(s^*(h,s^+),s^+\right) 
  \frac{\partial  s^* (h,s^+)}{\partial h} \; ds^+.\nonumber
\end{align}
We note that in the above derivation, the first term in the second line vanishes since
$s^*(h,s^+)\to\infty$ in the limit $s^+\to s^0(h)$ and $J(s^-,s^+) $ is bounded.

%________________________________________________________________
\section{Derivation of main result 3} \label{sec:limit}
%________________________________________________________________

With $h$ distributed according to $f_H(h)$, given by \eqref{eq:p_h2}, we derive in this section an
asymptotic expression for $f_H(h)$ in the limit $h\to\infty$. Examining \eqref{eq:main_approx}, we note that
$\hat{h}(s^-,s^+)$ is large when $s^-$ and/or $s^+$ are small, and thus in the large $h$ limit one can
consider only the contributions of the terms proportional to $s_-^{-2}$ and $s_+^{-2}$,
\begin{align}
h \approx \frac{{\lambda^2}}{(s^-)^2} + \frac{{\lambda^2}}{(s^+)^2}.
\end{align}
In this case, we find
\begin{align}
s^0(h) &= \frac{\lambda}{\sqrt{h}},\\
s^*(h,s^+) &= \frac{\lambda s^+}{\left[(s^+)^2 h - {\lambda^2}\right]^{1/2}},\\
\frac{\partial}{\partial h}\left(s^*(h,s^+)\right) &= \frac{-\lambda(s^+)^3}{2\left[(s^+)^2 h - {\lambda^2}\right]^{3/2}} = \frac{-1}{2{\lambda^2}}[s^*(h,s^+)]^3.
\end{align}
Substituting these values into \eqref{eq:p_h2} and dropping the arguments for $s^*$, i.e.
$s^*(h,s^+)\mapsto s^*$, we find 
\begin{align}
 f_H(h) &=  - \int_{\sqrt{{\lambda^2}/h}}^{\infty} \left(
 \frac{3^7\left[p \rho(\lambda)\right]^5}{32 \pi^3} s^+ s^*(s^*+ s^+) e^{-\frac{[3p\rho(\lambda)]^2}{4\pi}\left[(s^*)^2 + (s^+)^2 +s^*s^+\right]}\right) \left(\frac{-(s^*)^3}{2{\lambda^2}}\right)ds^+\nonumber\\
  &=  \frac{3^7\left[p \rho(\lambda)\right]^5}{32 \pi^3} \frac{1}{2{\lambda^2}}\int_{\sqrt{{\lambda^2}/h}}^{\infty} \left( \left(s^+ (s^*)^5 + (s^+)^2 (s^*)^4\right) e^{-\frac{[3p\rho(\lambda)]^2}{4\pi}\left[(s^*)^2 + (s^+)^2 +s^*s^+\right]}\right)  ds^+.\nonumber
\end{align}
The change of variables $u =  (s^+)^2 h - {\lambda^2} $ transforms this into
\begin{align}\label{eq:scale}
f_H(h)  & =  \frac{3^7\left[p \rho(\lambda)\right]^5}{32 \pi^3}\frac{{\lambda^2}}{4} h^{-7/2} I(h),
\end{align}
where we have defined 
\begin{equation}
I(h) = \int_{0}^{\infty} \left( 1 + \frac{{\lambda^2}}{u}\right)^{5/2} (u^{1/2}+ {\lambda}  ) e^{-\varphi(u)/ {h}}du,\label{eq:expre1}
\end{equation}
and
\begin{align}\label{eq:f}
\varphi(u) = \frac{[3p\rho(\lambda)]^2}{4\pi}(u+{\lambda^2})\left(1 + \frac{\lambda}{\sqrt{u}} + \frac{\lambda^2}{u}\right). 
\end{align}
The distribution $f_H(h)$ in \eqref{eq:scale}  depends on $h$ through the power law $h^{-7/2}$ as well as $I(h)$. In Appendix \ref{sec:scale_appendix}, we show  that \eqref{eq:expre1} has the large-$h$ scaling $I(h)=\mathcal{O}\left(\frac{h^{3/2}}{p^3}\right)$. Combining this with \eqref{eq:scale}, we find $f_H(h) = \mathcal{O}\left(\frac{p^2}{h^2}\right)$ for large $h$.

%________________________________________________________________
\section{Large-$h$ scaling of  $I(h)$} \label{sec:scale_appendix}
%________________________________________________________________

We now study how $I(h)$ given by \eqref{eq:expre1} scales in the limit of large $h$. Recall that the limit of large $h$ corresponds to when an eigenvalue $\lambda_i$ has a nearest-neighboring eigenvalue that is very close (i.e., $|\lambda_{i}-\lambda_{i\pm j}|\ll 1$), which results in large values of $h_i$ and subsequently the error of the empirical eigenvector (i.e., large
$\|{\bu}_i-\tu_i\|\approx h_i/n$).

Out strategy for evaluating \eqref{eq:expre1} is to split the integral into three regions of integration, which are chosen based on studying the function $\varphi(u)$. Examining \eqref{eq:f} for limiting values of $u$, we find that the function $\varphi(u)$ approaches $+\infty$, both as $u \to 0$ and as $u\to \infty$, and has the minimum 
\begin{equation}
    \min_{u\in[0\infty)}\varphi(u) = \varphi(\lambda^2)  = \frac{27}{2\pi}[\lambda p\rho(\lambda)]^2,\label{eq:minny}
\end{equation}
which occurs at $u = {\lambda^2}$. For large $h$, there are two values of $u$ such that $\varphi(u) = h$. We refer to these values as $u_1(h)$ and $u_2(h)$, with $u_1(h) < u_2(h)$. Considering the limits $u\to0$ and $u\to\infty$, we find the asymptotic approximations
\begin{align}
   u_1(h) &\to \frac{[3p\rho(\lambda)]^2}{4\pi} \lambda^4h^{-1}, \label{eq:u1}\\
   u_2(h) &\to  \frac{4\pi}{[3p\rho(\lambda)]^2}  h \label{eq:u2}.
\end{align}
We will evaluate \eqref{eq:expre1} by dividing the integration into
three ranges, $I(h) = I_1(h) + I_2(h) + I_3(h)$, where we define
\begin{align}
I_1(h) &= \int_{0}^{u_1(h)} \left( 1 + \frac{{\lambda^2}}{u}\right)^{5/2} (u^{1/2}+ {\lambda} ) e^{-\varphi(u)/ {h}}du,\label{eq:I1}\\
I_2(h) &= \int_{u_1(h)}^{u_2(h)} \left( 1 + \frac{{\lambda^2}}{u}\right)^{5/2} (u^{1/2}+ {\lambda} ) e^{-\varphi(u)/ {h}}du,\label{eq:I2}\\
I_3(h) &= \int_{u_2(h)}^{\infty} \left( 1 + \frac{{\lambda^2}}{u}\right)^{5/2} (u^{1/2}+ {\lambda} ) e^{-\varphi(u)/ {h}}du.\label{eq:I3}
\end{align}

We now study the $h\to\infty$ scaling for integrals $I_1(h)$, $I_2(h)$, and $I_3(h)$. Beginning with \eqref{eq:I1}, we first note that for the range $u\in(0,u_1(h)]$ that 
\begin{align}
(u + {\lambda^2})^{5/2}(u^{1/2} + {\lambda} ) \leq (u_1 + {\lambda^2})^{5/2}(u_1^{1/2} + {\lambda} ).
\end{align}
It follows that
\begin{equation}
   \left( 1 + \frac{{\lambda^2}}{u}\right)^{5/2} (u^{1/2}+ {\lambda} ) \leq E u^{-5/2},
\end{equation}
where 
\begin{equation}
E(\lambda)  = (u_1(h) + {\lambda^2})^{5/2}(u_1(h)^{1/2} + {\lambda}).
\end{equation}
Note that $E(\lambda)  \approx \lambda^6$ as $h\to\infty$, since $u_1(h) \to0$. Similarly, since $u$ is positive, one finds
\begin{align}
   \varphi(u) &=  {\frac{[3p\rho(\lambda)]^2}{4\pi}}(u+{\lambda^2})[1 + ({\lambda^2}/u) + ({\lambda^2}/u)^{1/2}]  \nonumber\\
   &\geq  {\frac{[3p\rho(\lambda)]^2}{4\pi}}({\lambda^2})({\lambda^2}/u) \nonumber\\
   &= Fu^{-1},
\end{align}
where we have defined 
\begin{equation}
    F =  \frac{[3p\rho(\lambda)]^2}{4\pi}{\lambda^4}.
\end{equation}
Using these two inequalities, we have
\begin{align}
I_1(h) &\leq   {E(\lambda)} \int_{0}^{u_1(h)} u^{-5/2} e^{-F/( {h} u)}du,\\
&= {E(\lambda)}   { \left(\frac{h}{F}\right)^{3/2}}\int_{F/( {h} u_1(h))}^{\infty} w^{1/2} e^{-w}dw\label{eq:blah}
\end{align}
which uses the change of variables $w = F/({h u(h)})$. Using \eqref{eq:u1}, the lower limit of integration converges as  $F/( {h} u_1(h)) \to 1$ with $h\to\infty$. The integral in \eqref{eq:blah} therefore limits to a constant, implying that $I_1(h)$ is dominated by a term which scales like $h^{3/2}$.

To estimate $I_3(h)$, note for large $h$ that that \eqref{eq:u2} implies $u> \lambda^2$ for any $u>u_2(h)$. It follows that
\begin{align}
\left( 1 + \frac{{\lambda^2}}{u}\right)^{5/2} (u^{1/2}+ {\lambda}) \leq 2^{5/2}(2 u^{1/2}).
\end{align} 
The integral $I_3(h)$ is thus dominated by 
\begin{align}
I_3(h) &\leq  {8}\int_{u_2(h)}^{\infty} u^{1/2}e^{-\varphi(u)/ {h}}du, \\
 &\leq  {8}\int_{u_2(h)}^{\infty} u^{1/2}\exp\left({-{\frac{[3p\rho(\lambda)]^2}{4\pi}}\frac{u}{h}}\right)du, 
\end{align}
where the second inequality uses $u>0$ and $\lambda^2/u>0$ to bound
\begin{align}
\varphi(u) &= {\frac{[3p\rho(\lambda)]^2}{4\pi}}(u+{\lambda^2})[1 + ({\lambda^2}/u) + ({\lambda^2}/u)^{1/2}] \nonumber\\
&\geq {\frac{[3p\rho(\lambda)]^2}{4\pi}} u.
\end{align}
We define the change of variables $w = {\frac{[3p\rho(\lambda)]^2}{4\pi}}\frac{u}{h}$ to obtain
\begin{align}
   I_3(h) \leq 8 \left({\frac{[3p\rho(\lambda)]^2}{4\pi}}\right)^{3/2}  
   \int_{{\frac{[3p\rho(\lambda)]^2}{4\pi}} {u_2(h)}/{h}}^{\infty} w^{1/2}e^{-w} dw. \label{eq:baba}
\end{align}
From \eqref{eq:u2}, the lower limit of integration converges as ${\frac{[3p\rho(\lambda)]^2}{4\pi}} {u_2(h)}/{h} \to 1$ and the integral in \eqref{eq:baba} converges to a constant as $h \to \infty$. Therefore $I_3(h)$ is also bounded by a term scaling as $h^{3/2}$.

We will now show that $I_2(h)$ has scaling $\mathcal{O}(h^{3/2})$ (as opposed to the other terms, which we showed are bounded by terms that scale as $h^{3/2}$).  Note that because of our definition of $u_1$ and $u_2$, and using that $\varphi(u)$ reaches a minimum at $u = {\lambda^2}\in(u_1,u_2)$, we find the bounds
\begin{align}
   \varphi(\lambda^2)/h \le \varphi(u)/h \le 1
\end{align}
for any $u \in(u_1, u_2)$. Substituting these into \eqref{eq:I3}, we bound $I_2(h)$ as
\begin{align}\label{eq:expre3}
   Q( {h})  e^{-1} \leq I_2(h) \leq Q( {h}) e^{-\varphi({\lambda^2})/ {h}},
\end{align}
where we have defined
\begin{align}\label{eq:Q}
Q( {h}) \equiv \int_{u_1( {h})}^{u_2( {h})} \left( 1 + \frac{{\lambda^2}}{u}\right)^{5/2} (u^{1/2}+ {\lambda}  )du.
\end{align}
Using the asymptotic approximations for $u_1(h) $ and $u_2(h) $ given by \eqref{eq:u1} and \eqref{eq:u2}, we  integrate \eqref{eq:Q} using the software Mathematica
(using the ``Series[$Q(h)$, \{h, Infinity,1\}]'' command) to obtain its asymptotic behavior, 
\begin{align}
    Q(h) &\approx  \frac{2^4\pi^{3/2}}{3^4[p\rho(\lambda)]^{3}} h^{3/2}.
\end{align}
Furthermore, we combine $\varphi({\lambda^2})/ {h}\to 0 $ with \eqref{eq:expre3} to obtain the asymptotic bound
\begin{align}\label{eq:expre3_b}
    Q( {h})  e^{-1} \leq I_2(h) \leq Q(h) .
\end{align}
We combine \eqref{eq:expre3_b} with \eqref{eq:blah} and \eqref{eq:baba} to obtain the large-$h$ scaling  $I(h)=\mathcal{O}\left(\frac{h^{3/2}}{p^3}\right)$.

\end{appendices}

\end{document}